\documentclass[amssymb, 11pt]{amsart}
\usepackage{latexsym}

\newcommand{\qb}[2]{{\left [{#1 \atop #2} \right]}}
\newlength{\standardunitlength}
\newtheorem{prop}{Proposition}[section]

\newtheorem{lemma}[prop]{Lemma}
\newtheorem{cor}[prop]{Corollary}
\newtheorem{theorem}[prop]{Theorem}

\newcommand{\GL}{\mathrm{GL}}
\newcommand{\tr}{\tau}
\newcommand{\diag}{\mathrm{diag}}
\newcommand{\Sp}{\mathrm{Sp}}
\newcommand{\SL}{\mathrm{SL}}
\newcommand{\PGL}{\mathrm{PGL}}

\begin{document}

\title [Conjugacy Classes] {Conjugacy class properties of the
extension of
$GL(n,q)$ generated by the inverse transpose involution}


\author{Jason Fulman}
\address{University of Pittsburgh\\
Pittsburgh, PA}
\email{fulman@math.pitt.edu}

\author{Robert Guralnick}
\address{University of Southern California \\ Los Angeles, CA}
\email{guralnic@math.usc.edu}



\date{March 17, 2003}

\thanks{Fulman was partially supported by National Security Agency
grant MDA904-03-1-0049.  Guralnick was partially supported by National
Science Foundation Agency grant DMS-0236185. }

\begin{abstract} Letting $\tau$ denote the inverse transpose automorphism of $GL(n,q)$,
a formula is obtained for the number of $g$ in $GL(n,q)$ so that
$gg^{\tau}$ is equal to a given element $h$. This generalizes a result
of Gow and Macdonald for the special case that $h$ is the identity. We
conclude that for $g$ random, $gg^{\tau}$ behaves like a hybrid of
symplectic and orthogonal groups. It is shown that our formula works
well with both cycle index generating functions and asymptotics, and
is related to the theory of random partitions. The derivation makes
use of models of representation theory of $GL(n,q)$ and of symmetric
function theory, including a new identity for Hall-Littlewood
polynomials. We obtain information about random elements of finite
symplectic groups in even characteristic, and explicit bounds for the
number of conjugacy classes and centralizer sizes in the extension of
$GL(n,q)$ generated by the inverse transpose automorphism.  We give a
second approach to these results using the theory of bilinear forms
over a field. The results in this paper are key tools in forthcoming
work of the authors on derangements in actions of almost simple
groups, and we give a few examples in this direction.
\end{abstract}

\maketitle

\section{Introduction} \label{intro}

 Let $F$ be a field and $G=GL(n,F)=GL(V)$.
 Let $g'$ denote the transpose of $g$ and let $g^{\tau}=(g')^{-1}$.
 Let $G^+ = \langle G, \tau \rangle$
 (so $\tau$ is an involution with $\tau g \tau = g^{\tr}$).
In the case of a finite field $F_q$ of size $q$, we will write $GL(n,q)$ or $G^+(n,q)$.

 The problem of counting the number of solutions to
 $gg^{\tau}=h$ where $g \in GL(n,q)$ and $h$ is a fixed element
 of $GL(n,q)$ has been addressed in several papers. It was
 proved by Gow \cite{gow1} in odd characteristic and later by Howlett
 and Zworestine \cite{HZ} in general that the number of
 solutions to this equation is equal to \[ \sum_{\chi \in
 Irr(GL(n,q))} \chi(h),\] where the sum is over all irreducible
 characters of $GL(n,q)$. For the special case when $h$ is
 the identity, it was proved by Gow \cite{gow1} in odd
 characteristic and later by Macdonald \cite{M} (pages 289-290) in
 general that this sum is equal to \[ (q-1)q^2(q^3-1)q^4(q^5-1)
 \cdots \] with $n$ factors altogether.

 One of the main results of this paper is a generalization of the
 formula of Gow and Macdonald to arbitrary elements $h$ in
 $GL(n,q)$. This question is of intrinsic interest since studying the
 conjugacy class statistics of $gg^{\tau}$ with $g$ random in
 $GL(n,q)$ is a natural cousin of studying conjugacy statistics of
 random elements in finite classical groups. The study of conjugacy
 classes of random elements in finite classical groups is a
 fascinating subject (see the survey \cite{F1}) and was crucial to our
 recent proof (see the series of papers beginning with \cite{FG1} and
 cited there) of a conjecture of Shalev stating that a finite simple
 group acting nontrivially on a finite set has at least a proportion
 of $\delta$ derangements (i.e. fixed point free elements) where
 $\delta>0$ is a universal constant (see also the paper \cite{B} of
 Boston et al asking a similar question).  The validity of Shalev's
 conjecture has applications to maps between varieties over finite
 fields and to random generation of groups.  It led us to investigate
 the proportion of derangements in almost simple groups (that is
 groups $H$ with $G \subseteq H \subseteq Aut(G)$ where $G$ is simple)
 and more particularly to the proportion of derangements in a given
 coset of the simple group.  One very special case of this set-up is
 when $H$ is $G^+(n,q)$ (or more precisely the quotient of $G^+$
 modulo scalars).  Then certainly if $g \tau$ is not a derangement in
 the action of $G^+(n,q)$ on a finite set $X$, then neither is $(g
 \tau)^2=gg^{\tau}$ (and if $gg^{\tau}$ has an odd number of fixed
 points, then $g \tau$ has a fixed point).  Hence understanding the
 behavior of $gg^{\tau}$ is important for the derangement
 problem. Another reason our enumeration is useful for the derangement
 problem is that in Section \ref{mincentsize} we obtain as a corollary
 a lower bound for the centralizer sizes of elements of
 $G^+(n,q)$. The sequel \cite{FG2} applies results in this article to
 the analog of Shalev's conjecture for almost simple groups,
 classifying (in a precise and quantitative way) how and when it
 fails. However here (Section \ref{derangements}) we at least give a
 few examples of how the tools in this paper can be used to study
 derangements, including examples where the proportion of derangements
 goes to $0$ as $q \rightarrow \infty$.

    In fact both for intrinsic interest and for applications to
 the derangement problem, it is useful to understand the asymptotics
 of conjugacy classes of $gg^{\tau}$ where $g$ is random in $GL(n,q)$.
 By this we mean the following.  Recall \cite{H} that the conjugacy
 classes of $GL(n,q)$ are parameterized by rational canonical form
 (i.e. an analog of Jordan form over finite fields): that is to each
 monic irreducible polynomial $\phi$ with coefficients in $F_q$, there
 is a association a partition $\lambda_{\phi}$ of size at most $n$ and
 conjugacy classes of $GL(n,q)$ corresponding to collections of
 partitions $\{\lambda_{\phi}\}$ satisfying the conditions that
 $|\lambda_z|=0$ and $\sum_{\phi} deg(\phi) |\lambda_{\phi}|=n$.  Here
 $|\lambda|$ denotes the size of a partition and $deg$ denotes the
 degree of the polynomial $\phi$.  We show that for any fixed finite
 collection of polynomials $S$, keeping $q$ fixed and letting $n
 \rightarrow \infty$, the partitions $\lambda_{\phi}(gg^{\tau})$ for
 $g$ random in $GL(n,q)$ are asymptotically independent for different
 polynomials in $S$, and we calculate their limit distributions.  We
 find (quite remarkably) that these limit distributions are
 essentially those defined and studied in \cite{F2} for the finite
 symplectic and orthogonal groups; this is one sense in which
 $gg^{\tau}$ behaves like a hybrid of symplectic and orthogonal
 groups. We also show that $gg^{\tau}$ with $g$ random in $GL(n,q)$
 has a cycle-index generating function. Both of these facts are
 crucial for asymptotic analysis.

 There are several ingredients in our method for evaluating the sum
\[ \sum_{\chi \in Irr(GL(n,q))} \chi(h).\]   First, we use work of
Klyachko \cite{K} (see also \cite{IS}) on models of irreducible
characters of $GL(n,q)$; that is a set of (not necessarily
irreducible) representations $\Theta_1,\cdots,\Theta_r$ of $GL(n,q)$
such that $\Theta_1+\cdots+\Theta_r$ is equivalent to the sum of all
irreducible representations of $GL(n,q)$, each occurring exactly
once. The difficult step in computing this sum in fact lies with
unipotent elements.  We solve this problem by translating it into the
language of Hall-Littlewood polynomials $P_{\lambda}(x;t)$, and then
establishing some new identities about these polynomials. For example
we prove that
\[ \sum_{\lambda} \frac{c_{\lambda}(t) P_{\lambda}(x;t)}{t^{o(\lambda)/2+|\lambda|/2}}
= \prod_{i \geq 1} \frac{1+x_i/t}{1+x_i} \prod_{i \leq j} \frac{1-x_ix_j}{1-x_ix_j/t},\]
where the sum is over all partitions of all natural numbers in which
the {\it even} parts occur with {\it even} multiplicity and all other
notation is defined in Section \ref{Hallpoly}.
We also give a simple combinatorial proof of Kawanaka's identity \cite{Ka}
\[ \sum_{\lambda}
\frac{c_{\lambda}(t) t^{o(\lambda)/2}
P_{\lambda}(x;t)}{t^{|\lambda|/2}} = \prod_{i \leq j}
\frac{1-x_ix_j}{1-x_ix_j/t},\] where the sum is over all partitions
of all natural numbers in which the {\it odd} parts occur with {\it even}
multiplicity. (Kawanaka's argument used Green's functions and work of
Lusztig on symmetric spaces).

 For the special case of unipotent elements $h$, there is
  another way (using representation theory of $GL(n,q)$ but
  nothing about models of irreducible representations)
to compute the number of $g$ such that $gg^{\tau}=h$.
This is a generalization of Macdonald's approach \cite{M}.
Our enumeration then follows from the new identity mentioned in
the previous paragraph. One nice aspect of this approach is that
it implies that when $h$ is unipotent, the formula for the number of
$g$ such that $gg^{\tau}=h$ is independent of whether the characteristic is
even or odd (we in fact use this observation in the argument of the previous paragraph).

    Section \ref{symplectic} studies a character sum \[ \sum_{\chi \in
    Irr(GL(2n,q)) \atop \chi \ even} \chi(h) \] where the sum is over
    a subset of all irreducible characters of $GL(n,q)$, defined more
    precisely in Section \ref{symplectic}. It follows from \cite{IS}
    or \cite{BKS} that this sum of irreducible characters is equal to
    the character obtained by inducing the trivial character of
    $Sp(2n,q)$ to $GL(2n,q)$. This induced character essentially tells
    us the proportion of elements of $Sp(2n,q)$ with a given rational
    canonical form (as an element of $GL(2n,q)$). Formulas for this
    proportion also follow from work of Wall \cite{W} on sizes of
    conjugacy classes in symplectic groups (see \cite{Ka} or \cite{F2}
    for a discussion of odd characteristic). However it is not obvious
    from Wall's treatment that when $h$ is unipotent, this proportion
    has the same form (as a function of $q$) in odd and even
    characteristic. A main result of Section \ref{symplectic} is a
    proof of this fact.

    The second part of this paper shifts to the viewpoint of
linear algebra.  Using bilinear forms, we give another approach to
enumerating $g$ so that $gg^{\tau}$ is equal to a given $h$.  While
this approach is not easy to work with in all cases (for instance if
$h$ is unipotent), it is more conceptual and (as with the
combinatorial approach) leads in all cases to lower bounds on
centralizer sizes of elements of $G^+(n,q)$ in the coset $GL(n,q)
\tau$. It is also quite convenient for treating a variation for $SL$.
We give explicit and useful (\cite{FG2}) upper bounds on the number of
conjugacy classes in $G^+(n,q)$ and also for the split extension of
$SL(n,q)$ generated by $\tau$.

 The precise organization of this paper is as follows.  Section
 \ref{Hallpoly} gives background on Hall-Littlewood polynomials and
 proves a number of identities about them. We obtain a new identity
 and also give an entirely combinatorial proof of an identity of
 Kawanaka on symmetric functions, avoiding work of Lusztig and Green's
 functions.  Section \ref{unipotent} generalizes Macdonald's approach
 to enumerating $g$ such that $gg^{\tau}=h$ in the case that $h$ is
 unipotent; for this purpose the identities of Section \ref{Hallpoly}
 are crucial. Section \ref{proof1} enumerates for arbitrary $h$, the
 number of $g$ satisfying $gg^{\tau}=h$.  It does this by using models
 of irreducible representations of $GL(n,q)$, converting the problem
 to one about Hall-Littlewood polynomials and using results from
 Section \ref{Hallpoly}. We emphasize that Section \ref{proof1} is
 independent of Section \ref{unipotent} in odd characteristic, and
 that the only fact used from Section \ref{unipotent} is that the
 enumeration for $h$ unipotent is independent (as a function of $q$)
 of whether the characteristic is even or odd.  Section
 \ref{symplectic} applies the same circle of ideas to the study of
 unipotent conjugacy classes in finite symplectic groups.  Section
 \ref{mincentsize} focuses on another corollary (very useful for the
 derangement problem), namely a lower bound on the centralizer size of
 an element of $G^+(n,q)$. Section \ref{cycle} shows that the
 enumeration of Section \ref{proof1} works well with both cycle index
 generating functions and asymptotics, and gives connections with the
 theory of random partitions. Section \ref{bilinearforms} consider the
 enumeration of $g$ such that $gg^{\tau}=h$ and its corollaries from
 the viewpoint of bilinear forms.  This gives an alternative proof of
 the enumeration for some $h$ and gives a different approach to lower
 bounds on centralizer sizes in Section \ref{mincentsize}.  Section
 \ref{numberclasses} provides an explicit upper bound on the number of
 $G^+(n,q)$ conjugacy classes in the coset $GL(n,q) \tau$ and also for
 the number of $SL(n,q)$ classes. Section \ref{derangements} gives a
 few examples of how tools in this paper can be used to study
 derangements in actions of $G^+(n,q)$.

\section{Identities for Hall-Littlewood Polynomials} \label{Hallpoly}

 To begin we collect some notation about partitions, much of it standard \cite{M}.
 Let $\lambda$ be a partition of some nonnegative integer
 $|\lambda|$ into parts $\lambda_1 \geq \lambda_2 \geq \cdots$.
 The symbol $m_i(\lambda)$ will denote the number of parts of $\lambda$
 of size $i$, and $\lambda'$ is the partition dual to $\lambda$ in the
 sense that $\lambda_i'=m_i+m_{i+1}+\cdots$. Let $n(\lambda)=\sum_i {\lambda_i' \choose 2}$.
 Let $l(\lambda)$ denote the number of parts of $\lambda$ and $o(\lambda)$
 the number of odd parts of $\lambda$.

 It is often helpful to view partitions diagrammatically.
 The diagram associated to $\lambda$ is the set of ordered pairs $(i,j)$ of
 integers such that $1 \leq j \leq \lambda_i$.
 We use the convention that the row index $i$ increases as one
 goes downward and the column index $j$ increases as one goes across.
 So the diagram of the partition $(5,4,4,1)$ is
 \[ \begin{array}{c c c c c}
             \framebox{}& \framebox{}& \framebox{}& \framebox{} & \framebox{}  \\
     \framebox{}& \framebox{} & \framebox{} & \framebox{} &\\
     \framebox{}& \framebox{} & \framebox{} & \framebox{}& \\
             \framebox{} &&&&
       \end{array} \] If a partition $\lambda$
       contains a partition $\mu$, then $\lambda-\mu$
       denotes the boxes in $\lambda$ which are not in $\mu$.
       One calls $\lambda-\mu$ a vertical strip if all of its boxes are in different rows.

Let $s$ denote some box in the diagram of the partition of $\lambda$.
Then $a_{\lambda}(s)$ (the arm of $s$) will denote the number of
boxes in the diagram of $\lambda$ in the same row as $s$ and to the east of $s$.
Similarly $l_{\lambda}(s)$ (the leg of $s$) will
denote the number of boxes in the diagram of $\lambda$
in the same column of $s$ and to the south of $s$.
When the partition $\lambda$ is clear from context
we sometimes omit the $\lambda$. Then one defines
\[ c_{\lambda}(t) = \prod_{s \in \lambda: a(s)=0,l(s) \ even} (1-t^{l(s)+1}) \]
where the product is over boxes $s$ in $\lambda$ with $a(s)=0$ and $l(s)$ even.

 This paper shall use the Hall-Littlewood polynomials
 $P_{\lambda}(x_1,x_2,\cdots;t)$. We often abbreviate this as
 $P_{\lambda}(x;t)$. They interpolate between Schur functions ($t=0$)
 and monomial symmetric functions ($t=1$). These are discussed
 thoroughly in Chapter 3 of \cite{M}. For the convenience of the
 reader we recall the definition of these polynomials and several
 properties of them which will be needed. Let $\lambda$ be a partition
 with $n$ parts (some of which may equal 0). Letting $v_{\lambda}(t) =
 \prod_{i \geq 0} \prod_{j=1}^{m_i(\lambda)} \frac{1-t^j}{1-t} $,
 define \[ P_{\lambda}(x_1,\cdots,x_n;t) = \frac{1}{v_{\lambda}(t)}
 \sum_{w \in S_n} w \left( x_1^{\lambda_1} \cdots x_n^{\lambda_n}
 \prod_{i<j} \frac{x_i-tx_j}{x_i-x_j} \right).\] The Hall-Littlewood
 polynomial is what one obtains by letting the number of variables go
 to infinity. We also recall that any symmetric function has an
 expansion in terms of the Hall-Littlewood polynomials.

 We shall also need the notion of a Hall polynomial. The Hall
 polynomial $g^{\lambda}_{\mu,\nu}(p)$ is the number of
 subgroups $H$ of an abelian p-group $G$ of type $\lambda$ such
 that $H$ has type $\mu$ and $G/H$ has type $\nu$. This is a
 polynomial in $p$ when $\lambda,\mu,\nu$ are fixed. For
 further discussion, see Chapter 2 of \cite{M}.
 The elementary symmetric function $e_r(x) =
 \prod_{i_1<\cdots<i_r} x_{i_1} \cdots x_{i_r}$ will be used. The notation
 $\qb{n}{m}$ denotes the $q$-binomial coefficient $\frac{(q^n-1)
 \cdots (q-1)}{(q^m-1) \cdots (q-1) (q^{n-m}-1) \cdots (q-1)}$.

 The following facts about Hall-Littlewood polynomials are needed.
 We emphasize that the proofs of these lemmas are entirely combinatorial.

\begin{lemma} \label{HallisHall} (\cite{M}, Section 3.3)
\[ t^{n(\mu)} P_{\mu}(x;t) t^{n(\nu)} P_{\nu}(x;t) =
\sum_{\lambda} g_{\mu,\nu}^{\lambda}(1/t) t^{n(\lambda)} P_{\lambda}(x;t)\]
where the sum is over all partitions $\lambda$.
\end{lemma}

\begin{lemma} \label{Macident} (\cite{M}, page 219)
\[ \sum_{\lambda} t^{n(\lambda)} \prod_{j=1}^{l(\lambda)} (1+t^{1-j}y) P_{\lambda}(x;t)
= \prod_{j \geq 1} \frac{1+x_j y}{1-x_j}\] where the sum is over all partitions $\lambda$.
\end{lemma}

\begin{lemma} \label{MacHallSum} (\cite{M}, page 231)
\[ \sum_{\mu} P_{\mu}(x;t) = \prod_i \frac{1}{1-x_i^2} \prod_{i<j} \frac{1-tx_ix_j}{1-x_ix_j}\]
where the sum is over partitions $\mu$ with all parts {\it even}.
\end{lemma}

 We shall also employ the following Pieri type formula which says how to
 multiply a Hall-Littlewood polynomial by an elementary symmetric function.

\begin{lemma} \label{Pieri} (\cite{M}, page 341)
\begin{eqnarray*} & & P_{\mu}(x;t) e_r(x)\\ & =
& \sum_{\lambda} P_{\lambda}(x;t) \prod_{j \geq 1} \frac{(t^{\lambda_j'-\lambda_{j+1}'}-1)
\cdots (t-1)}{(t^{\lambda_j'-\mu_j'}-1)\cdots (t-1) (t^{\mu_j'-\lambda_{j+1}'}-1)
\cdots (t-1)} \end{eqnarray*}
where the sum is over $\lambda$ such that $\lambda-\mu$ is a vertical strip of size $r$.
\end{lemma}

 Next we need some $q$-series identities (Lemma \ref{qseries1}
 and \ref{qseries2}). For this we recall a result of Euler.

\begin{lemma} \label{Euler} (Euler, page 19 of \cite{A})
\[ 1+\sum_{j=1}^{\infty} \frac{u^j}{(1-1/q) \cdots (1-1/q^j)} =
\prod_{j=0}^{\infty} \frac{1}{1-u/q^j}.\]
\end{lemma}

\begin{lemma} \label{qseries1} Let $(1/q)_a$ denote $(1-1/q) \cdots (1-1/q^a)$.
Then the expression $\sum_{r=0}^n \frac{(-1)^{n-r} (1/q)_n q^r}{(1/q)_r (1/q)_{n-r}}$
is equal to
\[ \left\{ \begin{array}{ll} q^n (1-1/q)(1-1/q^3) \cdots (1-1/q^{n-1}) &
\mbox{if \ $n$ \ even}\\
q^n (1-1/q)(1-1/q^3) \cdots (1-1/q^n) & \mbox{if \ $n$ \ odd}
\end{array} \right. \]
\end{lemma}

\begin{proof} We consider a generating function for a slightly modified sum.
\begin{eqnarray*}
& & \sum_{n=0}^{\infty} u^n \sum_{r=0}^n \frac{(-1)^{n-r} q^r}{(1/q)_r (1/q)_{n-r}}\\
& = & \sum_{r=0}^{\infty} \frac{(-1)^r q^r}{(1/q)_r} \sum_{n=r}^{\infty} \frac{(-1)^n u^n}{(1/q)_{n-r}}\\
& = & \sum_{r=0}^{\infty} \frac{(-1)^r q^r}{(1/q)_r} \sum_{n=0}^{\infty} \frac{(-1)^{n+r} u^{n+r}}{(1/q)_{n}}\\
& = & \sum_{r=0}^{\infty} \frac{u^r q^r}{(1/q)_r} \sum_{n=0 }^{\infty} \frac{(-1)^n u^n}{(1/q)_{n}}\\
& = & \prod_{j=0}^{\infty} \frac{1}{1-uq/q^j} \prod_{j=0}^{\infty} \frac{1}{1+u/q^j}\\
& = & \frac{1}{1-uq} \prod_{j=0}^{\infty} \frac{1}{1-u^2/q^{2j}}\\
& = & \frac{1}{1-uq} \left(\sum_{j \geq 0} \frac{u^{2j}}{(1-1/q^2)(1-1/q^4) \cdots (1-1/q^{2j})} \right). \end{eqnarray*} The fourth and sixth equalities have used Lemma \ref{Euler}.

 The coefficient of $u^n$ in this generating function is
 \begin{eqnarray*}
& & q^n \sum_{s=0}^{\lfloor n/2 \rfloor} \frac{1}{q^{2s} (1-1/q^2) \cdots (1-1/q^{2s})}\\
& = & q^n \frac{1}{(1-1/q^2) \cdots (1-1/q^{2 \lfloor n/2 \rfloor})}
\end{eqnarray*}
where the equality is proved by induction.
Note that this establishes the lemma since the generating function
was for the sought sum divided by $(1/q)_n$. \end{proof}

\begin{lemma} \label{qseries2} (\cite{A}, page 37)
Let $(1/q)_a$ denote $(1-1/q) \cdots (1-1/q^a)$.
Then $\sum_{r=0}^n (-1)^r \frac{(1/q)_n}{(1/q)_r (1/q)_{n-r}}$ is equal to
\[ \left\{ \begin{array}{ll} (1-1/q)(1-1/q^3) \cdots
(1-1/q^{n-1}) & \mbox{if \ $n$ \ even}\\ 0 & \mbox{if \ $n$ \ odd}
\end{array} \right. \]
\end{lemma}

 Now we establish an identity for Hall-Littlewood polynomials which is at
 the heart of this paper.

\begin{theorem} \label{newHall}
\[ \sum_{\lambda} \frac{c_{\lambda}(t) P_{\lambda}(x;t)}{t^{o(\lambda)/2+|\lambda|/2}} =
\prod_{i
\geq 1} \frac{1+x_i/t}{1+x_i} \prod_{i
\leq j} \frac{1-x_ix_j}{1-x_ix_j/t}\] where the sum is over partitions $\lambda$
in which all ${\it even}$ parts occur with ${\it even}$ multiplicity.
\end{theorem}

\begin{proof} Throughout we replace $t$ by $1/q$.
Write the right-hand side as
\[ \prod_{i} (1+x_iq)(1-x_i) \prod_i \frac{1}{1-qx_i^2}
\prod_{i<j} \frac{1-x_ix_j}{1-qx_ix_j }.\]
By Lemma \ref{MacHallSum}, this is
\[ \prod_{i} (1+x_i q)(1-x_i) \sum_{\mu} q^{|\mu|/2} P_{\mu}(x;1/q)\]
where the sum is over $\mu$ with all parts even.

 Next let us consider the coefficient of $P_{\tau}(x;1/q)$ in
 \[ \prod_{i} (1-x_i) \sum_{\mu} q^{|\mu|/2} P_{\mu}(x;1/q)\] where
 the sum is over $\mu$ with all parts even.
 Note that $\prod_{i} (1-x_i)= \sum_{r \geq 0} (-1)^r e_r(x)$
 where the $e_r(x)$ are the elementary symmetric functions.
 The Pieri-type rule (Lemma \ref{Pieri}) says that the effect of
 multiplying by $e_r$ is to add a size $r$ vertical strip with weights
 depending on the vertical strip. Observe that from $\tau$ there is a
 unique way of removing a vertical strip so as to get a partition with all
 parts even-one simply reduces the odd parts by 1. Hence the coefficient of
 $P_{\tau}(x;1/q)$ in
 \[ \prod_{i} (1-x_i) \sum_{\mu} q^{|\mu|/2} P_{\mu}(x;1/q)\]
 (where the sum is over $\mu$ with all parts even) is equal to
 $(-1)^{o(\tau)} q^{\frac{|\tau|-o(\tau)}{2}}$.

 Thus we need to find the coefficient of
 $P_{\lambda}(x;1/q)$ in
 \[ \prod_{i} (1+x_i q) \sum_{\tau} (-1)^{o(\tau)}
 q^{\frac{|\tau|-o(\tau)}{2}} P_{\tau}(x;1/q) \]
 where the sum is over all partitions $\tau$.
 Since $\prod_i (1+x_iq) = \sum_{r \geq 0} q^r e_r(x)$,
 we can again use the Pieri-type rule (Lemma \ref{Pieri}).
 Here however there are many possible ways of removing vertical strips from
 $\lambda$ since there are no restrictions on $\tau$.
 In fact using the notation that $(1/q)_a=(1-1/q) \cdots (1-1/q^a)$,
 one sees by Lemma \ref{Pieri} that the sought coefficient of $P_{\lambda}$ is precisely

\begin{eqnarray*}
& & \prod_{j \ odd} \sum_{r=0}^{m_j(\lambda)} q^r \frac{(1/q)_{m_j(\lambda)}}{(1/q)_r (1/q)_{m_j(\lambda)-r}} (-1)^{m_j(\lambda)-r} q^{\frac{[j m_j(\lambda)-r-(m_j(\lambda)-r)]}{2}}\\
& & \cdot \prod_{j \ even} \sum_{r=0}^{m_j(\lambda)} q^r \frac{(1/q)_{m_j(\lambda)}}{(1/q)_r (1/q)_{m_j(\lambda)-r}} (-1)^{r} q^{\frac{[j m_j(\lambda)-r-(r)]}{2}}\\
& = & \prod_{j \ odd} q^{(j-1)m_j(\lambda)/2} \sum_{r=0}^{m_j(\lambda)} q^r \frac{(1/q)_{m_j(\lambda)}}{(1/q)_r (1/q)_{m_j(\lambda)-r}} (-1)^{m_j(\lambda)-r}\\
& & \cdot \prod_{j \ even} q^{j m_j(\lambda)/2} \sum_{r=0}^{m_j(\lambda)} \frac{(1/q)_{m_j(\lambda)}}{(1/q)_r (1/q)_{m_j(\lambda)-r}} (-1)^{r}. \end{eqnarray*}

 By Lemma \ref{qseries2}, this vanishes if some even part of $\lambda$
 has odd multiplicity. Otherwise by Lemmas \ref{qseries1} and \ref{qseries2} is it equal to \begin{eqnarray*}
& = & \prod_{j \ odd \atop m_j(\lambda) \ even} q^{(j-1)m_j(\lambda)/2} q^{m_j(\lambda)} (1-1/q) (1-1/q^3) \cdots (1-1/q^{m_j(\lambda)-1})\\
& & \cdot \prod_{j \ odd \atop m_j(\lambda) \ odd} q^{(j-1)m_j(\lambda)/2} q^{m_j(\lambda)} (1-1/q) (1-1/q^3) \cdots (1-1/q^{m_j(\lambda)})\\
& & \cdot \prod_{j \ even \atop m_j(\lambda) \ even} q^{jm_j(\lambda)/2} (1-1/q) (1-1/q^3) \cdots (1-1/q^{m_j(\lambda)-1})\\
& = & c_{\lambda}(1/q) q^{o(\lambda)/2+|\lambda|/2}, \end{eqnarray*} as desired. \end{proof}

 The following identity of Kawanaka (\cite{Ka}) will also be needed.
 As his proof used Green's functions and work of Lusztig on symmetric spaces,
 we give a combinatorial proof using the same method as the proof of Theorem \ref{newHall}.
 One lemma (essentially a reformulation of Lemma \ref{qseries1}) will be used.

\begin{lemma} \label{qseries3} Let $(1/q)_a$ denote $(1-1/q) \cdots (1-1/q^a)$.
Then the expression $\sum_{r=0}^n \frac{(-1)^{r} (1/q)_n }{q^r(1/q)_r (1/q)_{n-r}}$ is equal to
\[ \left\{ \begin{array}{ll} (1-1/q)(1-1/q^3) \cdots (1-1/q^{n-1}) & \mbox{if \ $n$ \ even}\\
(1-1/q)(1-1/q^3) \cdots (1-1/q^n) & \mbox{if \ $n$ \ odd} \end{array}
\right. \]
\end{lemma}

\begin{proof} Observe that
\begin{eqnarray*} & & \sum_{r=0}^n \frac{(-1)^{r} (1/q)_n }{q^r(1/q)_r (1/q)_{n-r}} \\
& = & \frac{1}{q^n} \sum_{r=0}^n \frac{(-1)^{r} q^{n-r} (1/q)_n}{(1/q)_r (1/q)_{n-r}}\\
& = & \frac{1}{q^n} \sum_{r=0}^n \frac{(-1)^{n-r} q^{r} (1/q)_n}{(1/q)_{n-r} (1/q)_{r}}.
\end{eqnarray*}
Now apply Lemma \ref{qseries1}.
\end{proof}

\begin{theorem} \label{Kalemma}  (\cite{Ka})
\[ \sum_{\lambda} t^{\frac{o(\lambda)-|\lambda|}{2}} c_{\lambda}(t) P_{\lambda}(x;t) =
\prod_{i \leq j} \frac{1-x_ix_j}{1-x_ix_j/t}\]
where the sum is over all partitions where the {\it odd} parts occur with
{\it even} multiplicity.
\end{theorem}

\begin{proof} Throughout we replace $t$ by $1/q$. Write the right-hand side as
\[ \prod_{i} (1+x_i)(1-x_i) \prod_i \frac{1}{1-q x_i^2} \prod_{i<j} \frac{1-x_ix_j}{1-qx_ix_j }.\]
>From the proof of Theorem \ref{newHall}, one sees that this is equal to
\[ \prod_{i} (1+x_i) \sum_{\tau} (-1)^{o(\tau)} q^{\frac{|\tau|-o(\tau)}{2}} P_{\tau}(x;1/q) \]
where the sum is over all partitions $\tau$. Since $\prod_i (1+x_i) = \sum_{r \geq 0} e_r(x)$,
we can use the Pieri-type rule (Lemma \ref{Pieri}).
What emerges is that the coefficient of $P_{\lambda}(x;t)$ is equal to
\begin{eqnarray*}
& & \prod_{j \ odd} \sum_{r=0}^{m_j(\lambda)} \frac{(1/q)_{m_j(\lambda)}}{(1/q)_r (1/q)_{m_j(\lambda)-r}} (-1)^{m_j(\lambda)-r} q^{\frac{[j m_j(\lambda)-r-(m_j(\lambda)-r)]}{2}}\\
& & \cdot \prod_{j \ even} \sum_{r=0}^{m_j(\lambda)} \frac{(1/q)_{m_j(\lambda)}}{(1/q)_r (1/q)_{m_j(\lambda)-r}} (-1)^{r} q^{\frac{[j m_j(\lambda)-r-(r)]}{2}}\\
& = & \prod_{j \ odd} q^{(j-1)m_j(\lambda)/2} \sum_{r=0}^{m_j(\lambda)} \frac{(1/q)_{m_j(\lambda)}}{(1/q)_r (1/q)_{m_j(\lambda)-r}} (-1)^{m_j(\lambda)-r}\\
& & \cdot \prod_{j \ even} q^{jm_j(\lambda)/2} \sum_{r=0}^{m_j(\lambda)} \frac{(1/q)_{m_j(\lambda)}}{q^r (1/q)_r (1/q)_{m_j(\lambda)-r}} (-1)^{r}. \end{eqnarray*}

 By Lemma \ref{qseries2}, this vanishes if some odd part of $\lambda$
 has odd multiplicity. Otherwise by Lemmas \ref{qseries2} and \ref{qseries3}
 is it equal to
 \begin{eqnarray*}
& = & \prod_{j \ even \atop m_j(\lambda) \ even} q^{jm_j(\lambda)/2} (1-1/q) (1-1/q^3) \cdots (1-1/q^{m_j(\lambda)-1})\\
& & \cdot \prod_{j \ even \atop m_j(\lambda) \ odd} q^{jm_j(\lambda)/2} (1-1/q) (1-1/q^3) \cdots (1-1/q^{m_j(\lambda)})\\
& & \cdot \prod_{j \ odd \atop m_j(\lambda) \ even} q^{(j-1)m_j(\lambda)/2} (1-1/q) (1-1/q^3) \cdots (1-1/q^{m_j(\lambda)-1})\\
& = & c_{\lambda}(1/q) q^{-o(\lambda)/2+|\lambda|/2}, \end{eqnarray*}
as desired.
\end{proof}

\section{Enumeration of $g$ such that $gg^{\tau}=h$ for $h$ Unipotent} \label{unipotent}

 This section finds a formula for the number of $g$
 such that $gg^{\tau}=h$ when $h$ is unipotent. This approach
 uses nothing about models of irreducible representations of
 $GL(n,q)$ and generalizes the approach used by Macdonald
 \cite{M} for the case when $h$ is the identity.
 We also use one of our identities about Hall-Littlewood polynomials from
 Section \ref{Hallpoly}. A different approach to the case of $h$ unipotent is
 given in Section \ref{proof1} (though we do use the fact from this section that the
 answer has the same form for odd and even characteristic).

 To proceed we require two lemmas.

\begin{lemma} \label{factorpoly} Let $N(q;d)$
denote the number of monic degree $d$ irreducible polynomials
with coefficients in $F_q$ and non-0 constant term. Then
\[ \prod_{d \geq 1} (1-u^d)^{-N(q;d)} = \frac{1-u}{1-uq}.\]
\end{lemma}

\begin{proof} Rewriting the sought equation as
\[ \frac{1}{1-u} \prod_{d \geq 1} (1-u^d)^{-N(q;d)} = \frac{1}{1-uq}\]
the result follows from unique factorization in the ring $F_q[x]$.
Indeed the coefficient of $u^n$ on the right hand side is $q^n$,
the total number of monic degree $n$ polynomials with coefficients in $F_q$.
The left hand side says that each such polynomial factors
uniquely into irreducible pieces.
\end{proof}

 In Lemma \ref{Schursum} $s_{\lambda}(x)$ denotes the Schur function.

\begin{lemma} \label{Schursum} (\cite{M}, page 76)
\[ \sum_{\lambda} s_{\lambda} = \prod_i \frac{1}{1-x_i} \prod_{i<j} \frac{1}{1-x_ix_j}\]
where the sum is over all partitions $\lambda$.
\end{lemma}

 Theorem \ref{repGL} uses the fact that the representation theory of
 $GL(n,q)$ can be understood entirely in terms of symmetric function theory.
 A full account of this can be found in Chapter 4 of \cite{M}.

\begin{theorem} \label{repGL} Let $h$ be a unipotent element in $GL(n,q)$
of type $\mu$ (thus its Jordan blocks have sizes equal to the part
sizes of $\mu$).  Then the proportion of $g$ in $GL(n,q)$ such that
$gg^{\tau}$ is conjugate to $h$ is $0$ unless all even parts of $\mu$
have even multiplicity. If all even parts of $\mu$ have even
multiplicity, then the proportion is
\[ \frac{1}{q^{n(\mu)+\frac{n}{2}-\frac{o(\mu)}{2}}
\prod_i (1-1/q^2) \cdots (1-1/q^{2 \lfloor m_i(\mu)/2 \rfloor})}.\]
\end{theorem}

\begin{proof} For this proof we assume familiarity with Chapter 4 of \cite{M} and
adhere to his notation. Thus $M_n$ denotes the nonzero elements of the
algebraic closure of $F_q$ which are fixed by the $nth$ power of the
Frobenius map and $L_n$ is the character group of $M_n$. Also $\Theta$
is the set of primitive $F$-orbits $\theta$ in $\cup_n L_n$ and
$deg(\theta)$ is the $n$ such that $\theta$ is a primitive orbit in
$L_n$. The irreducible representations of $GL(n,q)$ are parameterized
by all ways of associating partitions $\lambda(\theta)$ to each
element of $\Theta$ in such a way that $\sum_{\theta \in \Theta}
deg(\theta) |\lambda_{\theta}|=n$.

 Let $p_{\lambda}(x) = \prod_{r \geq 1} (\sum_i x_i^r)^{m_r(\lambda)}$
  be the $\lambda$ power sum symmetric function and let $z(\tau)$
  denote the centralizer size of an element of conjugacy class type
  $\tau$ in a symmetric group on $|\tau|$ symbols. Let
  $\omega^{\lambda}(\tau)$ denote the character of the symmetric group
  parameterized by $\lambda$ on the conjugacy class parameterized by
  $\tau$. Chapter 4 of \cite{M} implies that the character value of an
  irreducible representation of type $\{\lambda(\theta)\}$ on a
  unipotent element of type $\mu$ is $q^{n(\mu)}$ multiplied by the
  coefficient of $P_{\mu}(x;1/q)$ in the symmetric function
  \begin{eqnarray*} && \prod_{\theta \in \Theta} \sum_{\tau}
  \frac{1}{z_{\tau}} \omega^{\lambda(\theta)}(\tau) \prod_{r \geq 1}
  ((-1)^{r \cdot deg(\theta)-1} p_{r \cdot
  deg(\theta)}(x))^{m_r(\tau)}\\ & = & \prod_{\theta \in \Theta}
  \sum_{\tau} \frac{1}{z_{\tau}} \omega^{\lambda(\theta)}(\tau)
  (-1)^{l(\tau)} \prod_{r \geq 1}
  p_{r}((-x)^{deg(\theta)})^{m_r(\tau)}\\ & = & \prod_{\theta \in
  \Theta} \sum_{\tau} \frac{1}{z_{\tau}}
  \omega^{\lambda(\theta)}(\tau) (-1)^{l(\tau)+|\tau|} \prod_{r \geq
  1} p_{r}((-x)^{deg(\theta)+1})^{m_r(\tau)}.
\end{eqnarray*}
Observe that $(-1)^{l(\tau)+|\tau|}$ is the sign of a permutation with
conjugacy class corresponding to the partition $\tau$.
As $s_{\lambda}(x)=
\sum_{\tau} \frac{1}{z_{\tau}} \omega^{\lambda}(\tau) p_{\tau}(x)$
and tensoring the irreducible representation of $S_n$ corresponding to
$\lambda$ by the sign representation simply switches
$s_{\lambda}(x)$ to $s_{\lambda'}(x)$, the above expression simplifies to
\[ \prod_{\theta \in \Theta} s_{\lambda(\theta)'}
(-(-x_1)^{deg(\theta)},-(-x_2)^{deg(\theta)},\cdots).\]

 Thus the sum over all irreducible characters $\chi$ of $GL(n,q)$ of
 their values on $h$ is $q^{n(\mu)}$ multiplied by the coefficient of
 $P_{\mu}(x;1/q)$ in the symmetric function \[ \prod_{d \geq 1}
 \left(\sum_{\lambda} s_{\lambda'}
 (-(-x_1)^{deg(\theta)},-(-x_2)^{deg(\theta)},\cdots)\right)^{N(q;d)},
 \] where $N(q;d)$ denotes the number of irreducible degree $d$
 polynomials over the field $F_q$ with non-zero constant term (these
 are in bijection with degree $d$ elements $\theta$ of $\Theta$).
 Invoking Lemma \ref{Schursum}, and using the fact that summing over
 all $\lambda$ is the same as summing over all $\lambda'$, this
 simplifies to \begin{eqnarray*} & & \prod_{d} \left( \prod_{i}
 (1+(-x_i)^{d})^{-1} \prod_{i<j} (1-x_i^dx_j^d)^{-1}
 \right)^{N(q;d)}\\ & = & \prod_d \left( \prod_{i}
 \frac{1-(-x_i)^{d}}{1-x_i^{2d}} \prod_{i<j} (1-x_i^dx_j^d)^{-1}
 \right)^{N(q;d)}. \end{eqnarray*} Using Lemma \ref{factorpoly}, this
 becomes \begin{eqnarray*} & & \prod_i \frac{1-x_i^2}{1+x_i}
 \frac{1+x_iq}{1-x_i^2q} \prod_{i<j} \frac{1-x_ix_j}{1-x_ix_jq}\\ & =
 & \frac{1+x_iq}{1+x_i} \prod_{i \leq j} \frac{1-x_ix_j}{1-x_ix_jq}.
\end{eqnarray*}

 It follows from Theorem \ref{newHall} that $q^{n(\mu)}$ multiplied by
 the coefficient of $P_{\mu}(x;1/q)$ in this symmetric function is $0$
 unless all even parts of $\mu$ have even multiplicity and is \[
 c_{\mu}(1/q) q^{n(\mu)+|\mu|/2+o(\mu)/2} \] if all even parts of
 $\mu$ have even multiplicity. Hence this is precisely the number of
 $g$ such that $gg^{\tau}$ is equal to a given unipotent element of
 type $\mu$. To determine the proportion of $g$ (random in $GL(n,q))$
 such that $gg^{\tau}$ is unipotent of type $\mu$, one need only
 divide this by the $GL(n,q)$ centralizer size of a unipotent element
 of type $\mu$, which is known (see page 191 of \cite{M}) to be
 $q^{2n(\mu)+|\mu|} \prod_i (1/q)_{m_i(\mu)}$. The result follows.
 \end{proof}

 Corollary \ref{char2ok} is immediate from Theorem \ref{repGL}.

\begin{cor} \label{char2ok} Let $h$ be a unipotent element of type $\lambda$
(thus the Jordan block sizes are the parts of $\lambda$).
Then the number of $g$ such that $gg^{\tau}=h$, viewed as a function of $q$,
has the same form in odd and even characteristic. \end{cor}

\section{General Enumeration of $g$ such that $gg^{\tau}=h$} \label{proof1}

 The purpose of this section is to derive a formula for the number of
 $g$ in $GL(n,q)$ such that $gg^{\tau}=h$ where $h$ is a fixed element of $GL(n,q)$.
 As mentioned in the introduction, the number of such $g$ is equal to
 \[ \sum_{\chi \in Irr(GL(n,q))} \chi(h).\]
 In particular, viewed as a function of $h$ this number is constant on
 conjugacy classes.

In fact (as noted in \cite{gow1}) this number is 0 unless
$h$ is a real (i.e. conjugate to its inverse) element of $GL(n,q)$.
Indeed, $(gg^{\tau})^{-1}=g'g^{-1}$ and $(gg^{\tau})'=g^{-1}g'$.
Thus $(gg^{\tau})^{-1}$ and $(gg^{\tau})'$ are conjugate.
The result now follows since any element in $GL(n,q)$
(in particular $gg^{\tau}$) is conjugate to its transpose.

 To begin we translate the problem of counting $g$
 so that $gg^{\tau}=h$ into a problem about Hall-Littlewood polynomials.
The following result of Klyachko \cite{K} (see \cite{IS} for an algebraic proof)
simplifies our task. Given groups $H \subset G$ and a character
$\chi$ of $H$, the symbol $\chi_H^G$ will denote the induced character.
We also recall a product $\circ$ which allow one to take a character $u_1$ of $GL(k,q)$
together with a character $u_2$ of $GL(n-k,q)$ and get a character of $GL(n,q)$.
Let $P_{k,n-k}$ be the parabolic subgroup of $GL(n,q)$ consisting of elements $g$ equal to
$$\begin{pmatrix} g_{11} & g_{12} \\  0 & g_{22}  \end{pmatrix},$$
where $g_{11} \in GL(k,q)$ and
$g_{22} \in GL(n-k,q)$. Then $u(g) = u_1(g_{11}) u_2(g_{22})$ is a class function on $P_{k,n-k}$
and inducing it to $GL(n,q)$ gives a character of $GL(n,q)$, denoted by $u_1 \circ u_2$.

\begin{theorem} \label{Klyach} (\cite{K}) Let $\gamma_k$ be the
Gelfand-Graev character of $GL(k,q)$.  Let
$\sigma_{2l}=1_{Sp(2l,q)}^{GL(2l,q)}$ denote the character of
$GL(2l,q)$ obtained by inducing the trivial character from $Sp(2l,q)$.
Then $\sum_{k+2l=n} \gamma_k \circ \sigma_{2l}$ is equal to the sum
of all irreducible characters of $GL(n,q)$, each occurring exactly
once.
\end{theorem}

 To apply Theorem \ref{Klyach}, one needs to know three things: a
 formula for $\gamma_k$, a formula for $\sigma_{2l}$, and how to
 compute the product $\circ$ using Hall polynomials. Fortunately, all
 of this information is available.

 At this point we remind the reader the conjugacy classes of $GL(n,q)$
 are parameterized by sets of partitions $\{ \lambda_{\phi} \}$
 (one for each monic irreducible polynomial $\phi$) satisfying $|\lambda_{z}|=0$
 and $\sum_{\phi} deg(\phi) |\lambda_{\phi}|=n$.
 The conjugacy data for real elements satisfies further restrictions.
 Namely there is an involution on monic irreducible polynomials with non-zero
 constant term sending a polynomial $\phi$ to
 $\bar{\phi}=\frac{z^{deg(\phi)} \phi(z)}{\phi(0)}$.
 The $\phi$ invariant under this involution are called self-conjugate.
 Real elements are precisely those which satisfy the additional constraint
 that $\lambda_{\phi}=\lambda_{\bar{\phi}}$.

 For the remainder of this section, we use the notation:
\[ A(\phi,\lambda_{\phi},i) =
\left\{
\begin{array}{ll} |U(m_i(\lambda_{\phi}),q^{deg(\phi)/2})| & \ if \ \phi=\bar{\phi} \\
|GL(m_i(\lambda_{\phi}),q^{deg(\phi)})|^{1/2} & \ if \ \phi \neq \bar{\phi}
\end{array} \right.\]
We remind the reader that $|GL(n,q)|=q^{n^2}(1/q)_n$ and that the size of $U(n,q)$ is $(-1)^n |GL(n,-q)|$. We define
$B(\phi,\lambda_{\phi})$ as
\[ \left\{ \begin{array}{ll}
q^{deg(\phi)[\sum_{h<i} hm_h(\lambda_{\phi})m_i(\lambda_{\phi})+\frac{1}{2}
\sum_i (i-1)m_i(\lambda_{\phi})^2]} \prod_i A(\phi,\lambda_{\phi},i) & \ \phi \neq z \pm 1 \\
q^{n(\lambda_{z+1})+\frac{|\lambda_{z+1}|}{2}+\frac{o(\lambda_{z+1})}{2}}
\prod_i (1-\frac{1}{q^2}) \cdots (1-\frac{1}{q^{2 \lfloor
m_i(\lambda_{z+1})/2 \rfloor}}) & \ \phi = z+1 \\
q^{n(\lambda_{z-1})+\frac{|\lambda_{z-1}|}{2}-\frac{o(\lambda_{z-1})}{2}}
\prod_i (1-\frac{1}{q^2}) \cdots (1-\frac{1}{q^{2 \lfloor
m_i(\lambda_{z-1})/2 \rfloor}}) & \ \phi = z-1 \end{array} \right. \]
and
where $\lfloor x \rfloor$ is the largest integer not exceeding
$x$. {\it In characteristic 2 we use the convention that the
polynomial $z+1$ does not exist-one uses formulas for $z-1$ instead}.

\begin{theorem} \label{Hallcor}  For $g$ random in $GL(n,q)$,
the chance that $gg^{\tau}$ has rational canonical form data $\{ \lambda_{\phi} \}$
is $0$ unless
\begin{enumerate}
\item $\lambda_{\phi}=\lambda_{\bar{\phi}}$ for all $\phi$
\item All even parts of $\lambda_{z-1}$ have even multiplicity.
\item All odd parts of $\lambda_{z+1}$ have even multiplicity.
\end{enumerate} If these conditions hold, then the chance is
\[ \prod_{\phi} \frac{1}{B(\phi,\lambda_{\phi})}. \]
\end{theorem}

\begin{proof} Note that the first condition must hold since as
explained at the beginning of this section, $gg^{\tau}$ is real.
Suppose first that the characteristic is odd. We apply Theorem \ref{Klyach}.
The Gelfand-Graev character $\gamma_k$ is well known.
For a simple proof in the case of $GL(k,q)$ see \cite{HZ} where it is shown that if
$dim(fix(h_1))$ denotes the dimension of the fixed space of an element $h_1$ in $GL(k,q)$,
the Gelfand-Graev character of $GL(k,q)$ evaluated at $h_1$ is $0$ if $h_1$ is not unipotent
and is equal to \[(-1)^{k-dim(fix(h_1))}(q^{dim(fix(h_1))}-1) \cdots (q-1)\] if
$h_1$ is unipotent. In the case when $h_1$ is unipotent,
let $\mu$ denote the partition of $k$ equal to $\lambda_{z-1}(h_1)$.
It is straightforward to see that $dim(fix(h_1))$ is equal to $l(\mu)$,
the number of parts of $\mu$.

The value of $1_{Sp(2l,q)}^{GL(2l,q)}(h_2)$ is also known; by the
general formula for induced characters \cite{S} it is simply
$\frac{1}{|Sp(2l,q)|}$ multiplied by the number of elements in
$GL(2l,q)$ which conjugate $h_2$ to something in $Sp(2l,q)$. This in
turn is $\frac{|C_{GL(2l,q)}(h_2)|}{|Sp(2l,q)|}$ multiplied by the
number of elements in $Sp(2l,q)$ with rational canonical form equal to
that of $h_2$ (i.e. elements conjugate to $h_2$ in $GL(2l,q)$).  The
centralizer sizes in general linear groups are well known (see for
instance page 191 of \cite{M}): if an element $h_2$ has conjugacy data
$\{\lambda_{\phi}\}$, the centralizer size is
\[\prod_{\phi} q^{2n(\lambda_{\phi})+|\lambda_{\phi}|} \prod_i
(1/q)_{m_i(\lambda_{\phi})}.\] If $h_2$ is not real, no elements
of $Sp(2l,q)$ have the rational canonical form of $h_2$.
Otherwise, from formulas in \cite{W}, one sees (as in \cite{F2} or
\cite{Ka}) that the number of elements in $Sp(2l,q)$ with the same
rational canonical form as $h_2$ is \[ \frac{|Sp(2l,q)|
q^{-n(\lambda_{z-1})-\frac{|\lambda_{z-1}|}{2}-\frac{o(\lambda_{z-1})}{2}}}{
\prod_i (1-1/q^2)\cdots (1-1/q^{2 \lfloor
\frac{m_i(\lambda_{z-1})}{2} \rfloor}) \prod_{\phi \neq z - 1}
B(\phi,\lambda_{\phi})}.\]

Section 4.3 of \cite{M} explains how to compute the product $\circ$ using Hall polynomials.
Applying this to the expression $\sum_{k+2l=n} \gamma_k \circ \sigma_{2l}$ from Theorem \ref{Klyach},
it follows that the proportion of $g$ such that
$gg^{\tau}$ is conjugate to a (real) element $h$ with rational
canonical form data $\{ \lambda_{\phi} \}$ is equal
to
\begin{eqnarray*} && \frac{1}{q^{2n(\lambda_{z-1})+|\lambda_{z-1}|} \prod_i
(1/q)_{m_i(\lambda_{z-1})}} \sum_{k,l \atop k+2l=n} \sum_{|\mu|=k}
\sum_{|\nu|=2l} g_{\mu,\nu}^{\lambda_{z-1}}(q)\\  & & \cdot
(-1)^{k-l(\mu)} (q^{l(\mu)}-1) \cdots (q-1) \\ & & \cdot \frac{q^{2n(\nu)+|\nu|}
\prod_i (1/q)_{m_i(\nu)}} {q^{n(\nu)+\frac{|\nu|}{2}+\frac{o(\nu)}{2}}
\prod_i (1-1/q^2) \cdots (1-1/q^{2 \lfloor m_i(\nu)/2
\rfloor}) \prod_{\phi \neq z - 1} B(\phi,\lambda_{\phi})},
\end{eqnarray*}
where all odd parts of $\nu$ occur with even multiplicity. This in turn
is equal to
\begin{eqnarray*} && \frac{1}{q^{2n(\lambda_{z-1})+|\lambda_{z-1}|} \prod_i
(1/q)_{m_i(\lambda_{z-1})}}  \sum_{\mu}
\sum_{\nu} g_{\mu,\nu}^{\lambda_{z-1}}(q)\\ & & \cdot
(-1)^{|\mu|-l(\mu)}(q^{l(\mu)}-1) \cdots
(q-1) \frac{q^{n(\nu)+\frac{|\nu|}{2}-\frac{o(\nu)}{2}}
c_{\nu}(1/q)} {\prod_{\phi \neq z - 1}
B(\phi,\lambda_{\phi})}
\end{eqnarray*} where $c(\nu)$ is as in Section \ref{Hallpoly}
and the sum is over all partitions $\mu,\nu$ with the condition
that all odd parts of $\nu$ occur with even multiplicity.

 Applying Lemma \ref{HallisHall}, this is equal to the
coefficient of $P_{\lambda_{z-1}}(x;1/q)$ in
\begin{eqnarray*} &&
\frac{1}{q^{n(\lambda_{z-1})+|\lambda_{z-1}|} \prod_i
(1/q)_{m_i(\lambda_{z-1})}} \sum_{\mu}
\frac{P_{\mu}(x;1/q)}{q^{n(\mu)}} (-1)^{|\mu|-l(\mu)} \\ && \cdot (q^{l(\mu)}-1)
\cdots (q-1) \sum_{\nu} P_{\nu}(x;1/q)
\frac{q^{\frac{|\nu|}{2}-\frac{o(\nu)}{2}} c_{\nu}(1/q)} {\prod_{\phi
\neq z - 1} B(\phi,\lambda_{\phi})} \end{eqnarray*}
where all odd
parts of $\nu$ occur with even multiplicity. Applying Lemma
\ref{Macident} (with the substitutions $t=1/q$, $y=-q$, and replacing all $x_i$
by their negatives), this simplifies to the coefficient of $P_{\lambda_{z-1}}(x;1/q)$ in
\begin{eqnarray*} &&
\frac{1}{q^{n(\lambda_{z-1})+|\lambda_{z-1}|} \prod_i
(1/q)_{m_i(\lambda_{z-1})}} \prod_{i \geq 1} \frac{1+x_iq}{1+x_i}\\ & & \cdot
\sum_{\nu} P_{\nu}(x;1/q) \frac{q^{\frac{|\nu|}{2}-\frac{o(\nu)}{2}}
c_{\nu}(1/q)} {\prod_{\phi \neq z - 1}
B(\phi,\lambda_{\phi})}.
\end{eqnarray*}
Using Theorem \ref{Kalemma}
this reduces to
\[ \frac{1}{q^{n(\lambda_{z-1})+|\lambda_{z-1}|}
\prod_i (1/q)_{m_i(\lambda_{z-1})}}
\prod_{i \geq 1} \frac{1+x_iq}{1+x_i} \prod_{i \leq j} \frac{1-x_ix_j}{1-x_ix_jq} \frac{1}{\prod_{\phi
\neq z - 1} B(\phi,\lambda_{\phi})}.\]
The theorem now follows (for odd
characteristic) from Theorem \ref{newHall}.

 To deduce the result in even characteristic, begin as in
 odd characteristic. The only change in even characteristic is
 with the dependence of the formula on $\lambda_{z-1}$--more precisely,
 the formula for the number of elements in $Sp(2l,q)$ with the same rational
 canonical form as $h_2$ is not obviously given by the expression stated in
 the odd characteristic case (in fact as we shall see later in Section \ref{symplectic},
 the two formulas are the same-but as this is somewhat painful to
 see directly from \cite{W} we do not use it in this proof).
 However from Corollary \ref{char2ok} in Section \ref{unipotent}
 (which does not use this theorem in its proof), one sees that the
 number of $g$ such that $gg^{\tau}$ is equal to a unipotent element of type
 $\lambda$ in $GL(|\lambda|,q)$ depends on $q$ in a way independent of the characteristic.
 Hence Theorem \ref{Hallcor} is valid in even characteristic. \end{proof}

 Next we note some corollaries of Theorem \ref{Hallcor}.
 More consequences appear in Section \ref{mincentsize} and \ref{cycle}.

\begin{cor} \label{likeSpO}
\begin{enumerate}
\item Suppose that $n$ is even. Let $C$ be a conjugacy class of
$GL(n,q)$ with the property that $\lambda_{z-1}$ is empty (i.e. the
eigenvalue 1 does not occur). Then the chance that $gg^{\tau} \in C$
for $g$ random in $GL(n,q)$ is equal to the chance that a random
element of $Sp(n,q)$ has $GL(n,q)$ conjugacy
class equal to $C$.
\item Suppose that $n$ is odd. Let $C$ be a conjugacy class of
$GL(n,q)$ with the property that $|\lambda_{z-1}|=1$ (i.e. the
eigenvalue 1 occurs with multiplicity 1).
Then the chance that $gg^{\tau} \in C$ for $g$ random in $GL(n,q)$ is equal to the chance that a random
element of $Sp(n-1,q)$ has $GL(n-1,q)$ conjugacy class data
$\{\lambda_{\phi}\}$ equal to that of $C$ except for $\lambda_{z-1}$ which is made empty.
\end{enumerate}
\end{cor}

\begin{proof} Both parts follows from Theorem \ref{Hallcor} and
Wall's formulas \cite{W} for conjugacy class sizes in finite
symplectic groups.  (In fact part 1 of the corollary is essentially
true because of Klyachko's result (Theorem \ref{Klyach}) together with
the fact that the Gelfand-Graev character of $GL(n,q)$ vanishes off of
unipotent elements; the full power of the symmetric function
calculations used to prove Theorem \ref{Hallcor} is not
needed). \end{proof}

 Corollary \ref{regss} shows that the proportion of regular semisimple elements $gg^{\tau}$
 is equal to a corresponding proportion in the symplectic groups
 (which was studied in \cite{GL} and \cite{FNP}) and will be crucial for our
 work on derangements. See for instance the examples in Section \ref{derangements}.

\begin{cor} \label{regss}
\begin{enumerate}
\item The proportion of elements $g \in GL(2n,q)$ such that $gg^{\tau}$ is regular
semisimple is equal to the proportion of elements in $Sp(2n,q)$ which are regular semisimple.
\item The proportion of elements $g \in GL(2n+1,q)$ such that $gg^{\tau}$ is regular
semisimple is equal to the proportion of elements in $Sp(2n,q)$ which are regular semisimple.
\end{enumerate}
\end{cor}

\begin{proof} The element $gg^{\tau} \in GL(n,q)$ is regular
semisimple precisely when its characteristic polynomial is squarefree.
Note that since the element $gg^{\tau}$ is real, this implies that
if $n$ is even the eigenvalue 1 does not occur, and if $n$ is odd,
the eigenvalue 1 occurs with multiplicity 1. Moreover an element of a
symplectic group is regular semisimple precisely when its characteristic
polynomial is square free; this implies that the eigenvalue 1 does not occur.
Now use Corollary \ref{likeSpO}.
\end{proof}

\section{Character Sums and Unipotent Symplectic Elements} \label{symplectic}

    The main purpose of this section is to use character theory of
    $GL(n,q)$ to compute the proportion of elements of $Sp(2n,q)$
    which are unipotent and have given rational canonical form in
    $GL(2n,q)$. In the case of odd characteristic this can be (and has
    been) alternatively computed directly from formulas of Wall
    \cite{W} (see \cite{F2} or \cite{Ka}). We shall see that the
    formula which arises is independent of whether the characteristic is odd or even. As one
    corollary the results of \cite{F2}, \cite{F3} on random elements
    of finite symplectic groups are applicable in even characteristic
    as well. We shall also be able to write down an expression for the
    number of elements (not necessarily unipotent) of $Sp(2n,q)$ which
    have given rational canonical form data $\{ \lambda_{\phi} \}$.

        To begin, we recall a result of \cite{IS}, \cite{BKS}. We use
    the notation that $1_{Sp(2n,q)}^{GL(2n,q)}$ denotes the character
    of $GL(2n,q)$ obtained by inducing the trivial character of
    $Sp(2n,q)$. All other notation conforms to that of Section
    \ref{unipotent}. Note that the partitions $\lambda$ in our
    notation are dual to those in the notation of \cite{IS},
    \cite{BKS}.

\begin{theorem} \label{inglis} (\cite{BKS},\cite{IS}) $1_{Sp(2n,q)}^{GL(2n,q)}$
is equal to the sum over all irreducible characters of $GL(2n,q)$ which satisfy the
constraint that $\lambda(\theta)'$ has all parts even for all $\theta \in \Theta$. \end{theorem}

    As mentioned above, Theorem \ref{Spunipotent} is known in odd
    characteristic, by a very different method of proof.

\begin{theorem} \label{Spunipotent} The proportion of elements $h$ of $Sp(2n,q)$
which are unipotent and have $GL(2n,q)$ rational canonical form of type $\mu$ is $0$
unless all odd parts of $\mu$ occur with even multiplicity. If all odd parts of $\mu$ occur with even multiplicity, it is
\[ \frac{1}{q^{n(\mu)+n+\frac{o(\mu)}{2}}
\prod_i (1-1/q^2) \cdots (1-1/q^{2\lfloor \frac{m_i(\mu)}{2}
\rfloor})} .\] \end{theorem}

\begin{proof} As explained in the proof of Theorem \ref{Hallcor}, $1_{Sp(2n,q)}^{GL(2n,q)}(h)$
is equal to $C_{GL(2n,q)}(h)$ multiplied by the proportion of elements of $Sp(2n,q)$
with rational canonical form equal to that of $h$.

    Next we evaluate $1_{Sp(2n,q)}^{GL(2n,q)}$ by applying the
    technique of Theorem \ref{repGL} to the result of Theorem
    \ref{inglis}. We conclude that $1_{Sp(2n,q)}^{GL(2n,q)}(h)$ is
    equal to $q^{n(\mu)}$ multiplied by the coefficient of
    $P_{\mu}(x;1/q)$ in \[ \prod_{d \geq 1} \left(\sum_{\lambda \atop
    all \ parts \ even} s_{\lambda}
    (-(-x_1)^{deg(\theta)},-(-x_2)^{deg(\theta)},\cdots)\right)^{N(q;d)},
    \] where $N(q;d)$ denotes the number of irreducible degree $d$
    polynomials over the field $F_q$ with non-zero constant term.

     By the Schur function case ($t=0$) of Lemma \ref{MacHallSum} and then
    Lemma \ref{factorpoly}, this simplifies to $q^{n(\mu)}$ multiplied
    by the coefficient of $P_{\mu}(x;1/q)$ in \begin{eqnarray*} & &
    \prod_{d \geq 1} \left( \prod_{i \geq 1} (1-x_i^{2d})^{-1}
    \prod_{i<j} (1-x_i^dx_j^d)^{-1} \right)^{N(q;d)}\\ & = & \prod_{i
    \geq 1} \frac{1-x_i^2}{1-qx_i^2} \prod_{i<j}
    \frac{1-x_ix_j}{1-qx_ix_j} \\ & = & \prod_{i \leq j}
    \frac{1-x_ix_j}{1-qx_ix_j}. \end{eqnarray*} The result now follows
    from Lemma \ref{Kalemma}, and the well-known formula (already
    used several times in this paper) \[ |C_{GL(2n,q)}(h)| =
    q^{2n+2n(\mu)} \prod_i (1/q)_{m_i(\lambda)}.\] \end{proof}

    We next give a second proof of Theorem \ref{Spunipotent}. This
    proof uses Wall's work and earlier results in this paper, but not
    Theorem \ref{inglis}.

\begin{proof} (Second proof) In odd characteristic this follows from Wall's formulas.
Thus it is enough to
show that the formula for the number of unipotent elements $h$ with
given partition $\lambda_{z-1}$ is (as a function of $q$) independent
of the characteristic.  We prove this by induction on $n$. By
Corollary \ref{char2ok}, we know that the number of elements $g$ of
$GL(n,q)$ with $gg^{\tau}$ conjugate to $h$ is (as a function of $q$)
independent of the characteristic.  Looking back at the proof of
Theorem \ref{Hallcor}, one sees that the formula for the number of
elements $g$ with $gg^{\tau}$ conjugate to $h$ is a sum over pairs of
partitions $\mu,\nu$ where the $\nu$ term involves the number of
elements of $Sp(2|\nu|,q)$ which are unipotent of type $\nu$, and the
$\mu$ term is in a form independent of the characteristic. Moreover,
precisely one term in this sum corresponds to $|\nu|=2n$--namely when
$\nu=\lambda_{z-1}$, and all other terms involve $\nu$ of smaller size
so by induction have the same form as a function of $q$ in either odd
or even characteristic.  This proves the result. \end{proof}

    In general, we have the follow result, which is immediate from Theorem \ref{Spunipotent}
    and results of \cite{W}.

\begin{cor} \label{randomSp} Using the notation before the proof of Theorem
\ref{Hallcor} and the convention that in even characteristic the
polynomial $z+1$ does not exist, the chance (in either odd or even
characteristic) that $g$ a random element of $Sp(2n,q)$ has rational
canonical form data $\{ \lambda_{\phi} \}$ is $0$ unless
$\lambda_{\phi}=\lambda_{\bar{\phi}}$ for all $\phi$ and $\lambda_{z
\pm 1}$ have all odd parts occur with even multiplicity.  If
$\lambda_{\phi}=\lambda_{\bar{\phi}}$ for all $\phi$, and $\lambda_{z
\pm 1}$ have all odd parts occur with even multiplicity, the chance is
\[ \frac{1}{q^{n(\lambda_{z-1})+\frac{|\lambda_{z-1}|}{2}+\frac{o(\lambda_{z-1})}{2}}
\prod_i (1-1/q^2) \cdots (1-1/q^{2\lfloor \frac{m_i(\lambda_{z-1})}{2} \rfloor})}
\prod_{\phi \neq z - 1} \frac{1}{B(\phi,\lambda_{\phi})}.\]
\end{cor}

\section{Minimum Centralizer Sizes} \label{mincentsize}

 The purpose of this section is to obtain a lower bound on the
 centralizer size of elements of $G^+(n,q)$.
 In fact we restrict consideration to elements in the coset
 $GL(n,q) \tau$, since centralizer sizes of elements of $G^+(n,q)$ in $GL(n,q)$
 are at most double the $GL(n,q)$ centralizer size of elements in $GL(n,q)$
 (and the paper \cite{FG3} lower bounded the minimum centralizer sizes of elements
 in $GL(n,q)$). The bound in this section is crucial to our study of
 derangements in \cite{FG2}.

 The bound of this section is derived as a consequence of Theorem
 \ref{Hallcor}. In Section \ref{bilinearforms}, we give a different
 approach using the theory of bilinear forms.

 Two lemmas are required.

\begin{lemma} \label{NP0lem} (\cite{NP}) Suppose that $q \geq 2$.
Then
$\prod_{i \geq 1} (1-\frac{1}{q^i}) \geq 1-\frac{1}{q}-\frac{1}{q^2}$.
\end{lemma}

\begin{lemma} \label{centbounduni}
\begin{enumerate}
\item Let $\lambda$ be a partition in which all even parts have even
multiplicity. Then \[
q^{n(\lambda_{z-1})+\frac{|\lambda|}{2}-\frac{o(\lambda)}{2}} \prod_i
(1-1/q^2) \cdots (1-1/q^{2 \lfloor m_i(\lambda)/2 \rfloor}) \geq
q^{\lfloor |\lambda|/2 \rfloor} (1-1/q^2-1/q^4).\]
\item Let $\lambda$ be a partition in which all odd parts have even
multiplicity. Then \[
q^{n(\lambda_{z-1})+\frac{|\lambda|}{2}+\frac{o(\lambda)}{2}} \prod_i
(1-1/q^2) \cdots (1-1/q^{2 \lfloor m_i(\lambda)/2 \rfloor}) \geq
q^{ |\lambda|/2 } (1-1/q^2-1/q^4).\]
\end{enumerate}
\end{lemma}

\begin{proof} Consider the first assertion. From \cite{F2} it is known that
\[  \sum_{|\lambda|=n} \frac{1}{ q^{n(\lambda)+
\frac{|\lambda|}{2}-\frac{o(\lambda)}{2}} \prod_i (1-1/q^2) \cdots
(1-1/q^{2 \lfloor m_i(\lambda)/2 \rfloor})} \] (the sum is over all
partitions of size $n$ in which all even parts have even multiplicity)
is the coefficient of $u^n$ in $\frac{1+u}{\prod_{j \geq 1}
(1-u^2/q^{2j-1})}$. From Lemma \ref{Euler} this coefficient is \[
\frac{1}{q^{\lfloor n/2 \rfloor} (1-1/q^2) \cdots (1-1/q^{2 \lfloor
n/2 \rfloor})}.\] Hence any particular term in this sum is at most
\[ \frac{1}{q^{\lfloor n/2 \rfloor} (1-1/q^2) \cdots (1-1/q^{2 \lfloor n/2 \rfloor})}, \]
so the result follows from Lemma \ref{NP0lem}.

    The proof of the second assertion is similar, using the fact
    from \cite{F2} that \[ \sum_{|\lambda|=n} \frac{1}{
    q^{n(\lambda)+ \frac{|\lambda|}{2}+\frac{o(\lambda)}{2}}
    \prod_i (1-1/q^2) \cdots (1-1/q^{2 \lfloor m_i(\lambda)/2
    \rfloor})} \] (the sum is over all partitions of size $n$ in
    which all odd parts have even multiplicity) is the coefficient
    of $u^n$ in $\frac{1}{\prod_{j \geq 1}
    (1-u^2/q^{2j-1})}$. \end{proof}

 Now the main result of this section can be proved.

\begin{theorem} \label{centralizersize}
The $GL(n,q)$ centralizer size of an element in the coset
$GL(n,q) \tau$ is at least
$(1-1/q^2-1/q^4)^2 q^{\lfloor n/2 \rfloor} \left( \frac{1-1/q}{4e log_q(n)} \right)^{1/2}$.
\end{theorem}

\begin{proof} Let $g \tau$ be an element of $G^+(n,q)$
whose square is equal to $h \in GL(n,q)$.
Let $s(h)$ denote the number of elements of in the coset $GL(n,q) \tau$ whose square is $h$.
Note that if $z \in C_{GL(n,q)}(h)$, then $(z g \tau z^{-1})^2=h$. Thus
\begin{eqnarray*} s(h) & \geq & \frac{|C_{GL(n,q)}(h)|}
{|C_{GL(n,q)}(g \tau) \cap C_{GL(n,q)}(h)|}\\
& \geq & \frac{|C_{GL(n,q)}(h)|}{|C_{GL(n,q)}(g \tau)|}.
\end{eqnarray*}
Hence $C_{GL(n,q)}(g \tau)$ is at least the reciprocal of the
proportion of elements $x$ in $GL(n,q)$ such that $xx^{\tau}$ is conjugate to $h$.
Thus Theorem \ref{Hallcor} implies that
\[ C_{GL(n,q)}(g \tau) \geq \prod_{\phi} B(\phi,\lambda_{\phi}) \]
where $B(\phi,\lambda_{\phi})$ is defined before Theorem \ref{Hallcor}
and $\{ \lambda_{\phi} \}$ is the conjugacy class data of $h$.

 Let $2m$ be the degree of the part of the characteristic polynomial
 of $h$ which is relatively prime to $(z-1)^2$. It follows from [FG3]
 that \[ \prod_{\phi \neq z \pm 1} B(\phi,\lambda_{\phi}) \geq q^{m}
 \left( \frac{1-1/q}{4e log_q(2m)} \right)^{1/2}.\] Note that Lemma
 \ref{centbounduni} implies that \[ B(z-1,\lambda_{z-1})
 B(z+1,\lambda_{z+1})  \geq  (1-1/q^2-1/q^4)^2 q^{\lfloor
 |\lambda_{z-1}|/2 \rfloor + |\lambda_{z+1}|/2} .\] The result follows
 since $m+ \lfloor \frac{|\lambda_{z-1}|}{2} \rfloor + \frac{|\lambda_{z+1}|}{2} =
 \lfloor \frac{n}{2} \rfloor$. \end{proof}

\section{Generating Functions, Asymptotics, and Random Partitions} \label{cycle}

 This section consists of some important corollaries of Theorem
 \ref{Hallcor}. Corollary \ref{cycleindex} gives a cycle index
 generating function.  The cycle index is a very useful tool for
 studying properties of random matrices and sometimes allows one to
 obtain results out of reach by other methods; see the survey
 \cite{F1} or \cite{F3}.  It is used in our work on the derangement
 problem (see \cite{FG2} and also Section \ref{derangements} of this
 paper).

 In Corollary \ref{cycleindex}, the $x_{\phi,\lambda}$ are variables.
 Recall that $B(\phi,\lambda)$ was defined in Section
 \ref{proof1}. Although one can write down a single generating
 function, it is more useful for asymptotic purposes to treat
 separately the cases that $n$ is odd or even. Indeed, the size of the
 partition $\lambda_{z-1}(gg^{\tau})$ is equal to $n$ modulo 2.

\begin{cor} \label{cycleindex} Let $e=0$ if the characteristic is even and $e=1$
if the characteristic is odd. In the equations below, $\phi$ denotes a
    monic irreducible polynomial over $F_q$ the $\{
    \phi,\bar{\phi} \}$ denote conjugate (unordered) pairs of
    non-selfconjugate monic irreducible polynomials.

\begin{enumerate}

\item \begin{eqnarray*} & & 1+\sum_{n \geq 1}
\frac{u^{2n}}{|GL(2n,q)|} \sum_{g \in GL(2n,q)} \prod_{\phi}
x_{\phi,\lambda_{\phi}(gg^{\tau})}\\ & = & \left( \sum_{|\lambda| \
even \atop i \ even \Rightarrow m_i \ even} \frac{ x_{z-1,\lambda}
u^{|\lambda|}}{q^{n(\lambda)+\frac{|\lambda|}{2}-\frac{o(\lambda)}{2}}
\prod_i (1-1/q^2) \cdots (1-1/q^{2 \lfloor \frac{m_i(\lambda)}{2}
\rfloor})} \right)\\ & & \cdot \left( \sum_{\lambda \atop i \ odd
\Rightarrow m_i \ even} \frac{ x_{z+1,\lambda}
u^{|\lambda|}}{q^{n(\lambda)+\frac{|\lambda|}{2}+\frac{o(\lambda)}{2}}
\prod_i (1-1/q^2) \cdots (1-1/q^{2 \lfloor \frac{m_i(\lambda)}{2}
\rfloor})} \right)^{e}\\ & & \cdot \prod_{\phi = \bar{\phi} \atop \phi
\neq z \pm 1} \left( \sum_{\lambda} \frac{x_{\phi,\lambda}
u^{|\lambda| \cdot deg(\phi)}}{B(\phi,\lambda)} \right)
\prod_{\{\phi,\bar{\phi}\} \atop \phi \neq \bar{\phi}} \left( \sum_{\lambda}
\frac{x_{\phi,\lambda}x_{\bar{\phi},\lambda} u^{2|\lambda| \cdot
deg(\phi)}}{B(\phi,\lambda) B(\bar{\phi},\lambda)} \right)
\end{eqnarray*}

\item  \begin{eqnarray*} & & 1+\sum_{n \geq 0} \frac{u^{2n+1}}{|GL(2n+1,q)|}
\sum_{g \in GL(2n+1,q)} \prod_{\phi} x_{\phi,\lambda_{\phi}(gg^{\tau})}\\
& = & \left( \sum_{|\lambda| \ odd \atop i \ even \Rightarrow m_i \ even}
\frac{ x_{z-1,\lambda} u^{|\lambda|}}{q^{n(\lambda)+\frac{|\lambda|}{2}-\frac{o(\lambda)}{2}}
\prod_i (1-1/q^2) \cdots (1-1/q^{2 \lfloor
\frac{m_i(\lambda)}{2} \rfloor})} \right)\\
& & \cdot \left(
\sum_{\lambda \atop i \ odd \Rightarrow m_i \ even}
\frac{ x_{z+1,\lambda} u^{|\lambda|}}{q^{n(\lambda)+\frac{|\lambda|}{2}+\frac{o(\lambda)}{2}}
\prod_i (1-1/q^2) \cdots (1-1/q^{2 \lfloor
\frac{m_i(\lambda)}{2} \rfloor})} \right)^{e}\\ & & \cdot \prod_{\phi =
\bar{\phi} \atop \phi \neq z \pm 1} \left( \sum_{\lambda}
\frac{x_{\phi,\lambda} u^{|\lambda| \cdot deg(\phi)}}{B(\phi,\lambda)} \right)
\prod_{\{\phi,\bar{\phi}\} \atop \phi \neq \bar{\phi}} \left( \sum_{\lambda}
\frac{x_{\phi,\lambda}x_{\bar{\phi},\lambda} u^{2|\lambda| \cdot deg(\phi)}}{B(\phi,\lambda)
B(\bar{\phi},\lambda)} \right)
\end{eqnarray*}

\end{enumerate}
 \end{cor}

\begin{proof} Consider the first part. Note that $|\lambda_{z-1}(gg^{\tau})|$ must be even
(since $|\lambda_{z+1}|$ is even and all self-conjugate polynomials other than $z \pm 1$ have even degree).
The coefficient of $u^{2n} \prod_{\phi} x_{\phi,\lambda_{\phi}}$
on the left-hand side is the proportion of elements $g$ in $GL(2n,q)$
such that $gg^{\tau}$ has rational canonical form data
$\{\lambda_{\phi}\}$.  By Theorem \ref{Hallcor}, this is also the
coefficient of $u^{2n} \prod_{\phi} x_{\phi,\lambda_{\phi}}$ on the
right-hand side. The second assertion is similar.
\end{proof}

 Next we give some asymptotic consequences of the cycle index for
 the theory of random partitions. First we note that the theory of random
 partitions is quite interesting (see the surveys \cite{F1} or \cite{Ok}).
 The paper \cite{F2} defined a probability measure on the set of all partitions with the
 property that all odd parts occur with even multiplicity by the formula
 \[ M_{Sp,u}(\lambda) = \prod_{i=1}^{\infty} (1-\frac{u^2}{q^{2i-1}})
\frac{u^{|\lambda|}}{q^{n(\lambda)+\frac{|\lambda|}{2}+\frac{o(\lambda)}{2}}
\prod_i (1-1/q^2) \cdots (1-1/q^{2 \lfloor
\frac{m_i(\lambda)}{2} \rfloor})} \]
where $u$ is a parameter.
It also defined a probability measure on the set of all partitions with the property that
all even parts occur with even multiplicity by the formula
\[ M_{O,u}(\lambda) = \frac{\prod_{i=1}^{\infty}
(1-\frac{u^2}{q^{2i-1}})}{1+u}
\frac{u^{|\lambda|}}{q^{n(\lambda)+\frac{|\lambda|}{2}-\frac{o(\lambda)}{2}}
\prod_i (1-1/q^2) \cdots (1-1/q^{2 \lfloor \frac{m_i(\lambda)}{2}
\rfloor})} \] where $u$ is a parameter. Note that in both of these
definitions, the size of $\lambda$ is not fixed.

    In fact for asymptotic purposes, it is useful to refine the
measure $M_{O,u}$ into two measures $M_{O,u,even}$ and
$M_{O,u,odd}$. The measure $M_{O,u,even}$ is supported on all
partitions of even size in which all even parts occur with even
multiplicity and is defined there as \[ \prod_{i=1}^{\infty}
(1-\frac{u^2}{q^{2i-1}})
\frac{u^{|\lambda|}}{q^{n(\lambda)+\frac{|\lambda|}{2}-\frac{o(\lambda)}{2}}
\prod_i (1-1/q^2) \cdots (1-1/q^{2 \lfloor \frac{m_i(\lambda)}{2}
\rfloor})}. \] The measure $M_{O,u,odd}$ is supported on all
partitions of odd size in which all even parts occur with even
multiplicity and is defined there as \[ \prod_{i=1}^{\infty}
(1-\frac{u^2}{q^{2i-1}})
\frac{u^{|\lambda|-1}}{q^{n(\lambda)+\frac{|\lambda|}{2}-\frac{o(\lambda)}{2}}
\prod_i (1-1/q^2) \cdots (1-1/q^{2 \lfloor \frac{m_i(\lambda)}{2}
\rfloor})}. \]

 In fact, as we shall now see, these random partitions are related to
 the study of $gg^{\tau}$ where $g$ is random in $GL(n,q)$.  Fixing a
 polynomial $\phi$ and choosing $g$ random, the partition
 $\lambda_{\phi}(gg^{\tau})$ is a random partition. In fact with
 $\phi,q$ fixed and $n \rightarrow \infty$ (with the value of $n$
 modulo 2 specified) and $g$ random in $GL(n,q)$, the partition
 $\lambda_{\phi}(gg^{\tau})$ has a limit distribution which can be
 identified. Moreover, except for the fact that
 $\lambda_{\phi}=\lambda_{\bar{\phi}}$, these partitions will be
 asymptotically independent, which is useful for asymptotic
 calculations. To prove this, a combinatorial lemma is required.

\begin{lemma} \label{identity} Let $f=1$ in even characteristic and
$f=2$ in odd characteristic.
\begin{eqnarray*}
1-u^2 & = &
\left(\prod_{i \geq 1} (1-u^2/q^{2i-1}) \right)^f
\prod_{\phi=\bar{\phi} \atop \phi \neq z \pm 1} \prod_{i
\geq 1} (1+u^{deg(\phi)}/(-1)^i q^{i \cdot deg(\phi)/2}) \\ & & \cdot
\prod_{\{\phi,\bar{\phi}\} \atop \phi \neq \bar{\phi}} \prod_{i \geq 1}
(1-u^{2deg(\phi)}/q^{i \cdot deg(\phi)}),
\end{eqnarray*}
where the final product is over conjugate (unordered)
pairs of non selfconjugate monic irreducible polynomials.
\end{lemma}

\begin{proof} This equation is the reciprocal of the equation obtained by setting
all variables (other than $u$) equal to 1 in the index of the symplectic groups \cite{F3}. \end{proof}

 Now the main theorem can be stated. The use of auxiliary
 randomization (i.e. randomizing the variable $n$) is a mainstay of
 statistical mechanics known as the grand canonical ensemble. The
 second part of Theorem \ref{asymptotics} is an example of the
 principle of equivalence of ensembles: as $n$ gets large the system
 for fixed $n$ (microcanonical ensemble) behaves like the grand
 canonical ensemble. We say that an infinite collection of random
 varialbes is independent if any finite subcollection is.

\begin{theorem} \label{asymptotics}
\begin{enumerate}
\item Fix $u$ with $0<u<1$. Then choose a random even natural number $N$
with the probability of getting $2n$ equal to $(1-u^2)u^{2n}$. Choose $g$
uniformly at random in $GL(2n,q)$ and let $\Lambda_{\phi}(gg^{\tau})$ be
the partition corresponding to the polynomial $\phi$ in the rational
canonical form of $gg^{\tau}$. Then as $\phi$ varies, aside from the
fact that $\Lambda_{\phi}=\Lambda_{\bar{\phi}}$, these random
variables are independent with probability laws the same as those for
the symplectic groups in Theorem 1 of \cite{F2} except for the
polynomial $z-1$ which has the distribution $M_{O,u,even}$.

\item  Fix $u$ with $0<u<1$. Then choose a random odd natural number $N$
with the probability of getting $2n+1$ equal to $(1-u^2)u^{2n}$. Choose $g$
uniformly at random in $GL(2n+1,q)$ and let $\Lambda_{\phi}(gg^{\tau})$ be
the partition corresponding to the polynomial $\phi$ in the rational
canonical form of $gg^{\tau}$. Then as $\phi$ varies, aside from the
fact that $\Lambda_{\phi}=\Lambda_{\bar{\phi}}$, these random
variables are independent with probability laws the same as those for
the symplectic groups in Theorem 1 of \cite{F2} except for the
polynomial $z-1$ which has the distribution $M_{O,u,odd}$.

\item Choose $g$ uniformly at random in $GL(2n,q)$ and let
$\Lambda_{\phi}(gg^{\tau})$ be the partition corresponding to the
polynomial $\phi$ in the rational canonical form of $gg^{\tau}$. Let
$q$ be fixed and $n \rightarrow \infty$. Then as $\phi$ varies, aside
from the fact that $\Lambda_{\phi}=\Lambda_{\bar{\phi}}$, these random
variables are independent with probability laws the same as those for
the symplectic groups in Theorem 1 of \cite{F2} except for the
polynomial $z-1$ which has the distribution $M_{O,1,even}$.

\item Choose $g$ uniformly at random in $GL(2n+1,q)$ and let
$\Lambda_{\phi}(gg^{\tau})$ be the partition corresponding to the
polynomial $\phi$ in the rational canonical form of $gg^{\tau}$. Let
$q$ be fixed and $n \rightarrow \infty$. Then as $\phi$ varies, aside
from the fact that $\Lambda_{\phi}=\Lambda_{\bar{\phi}}$, these random
variables are independent with probability laws the same as those for
the symplectic groups in Theorem 1 of \cite{F2} except for the
polynomial $z-1$ which has the distribution $M_{O,1,odd}$.

\end{enumerate}
\end{theorem}

\begin{proof} The method of proof is analogous to that used for the classical groups
(see the survey \cite{F1}). We treat the case of $n$ even as the case of $n$ odd is similar.
For the first part, one multiplies the cycle index
(Corollary \ref{cycleindex}) by the equation in Lemma
 \ref{identity}. To prove the third assertion one uses the fact that
 if a Taylor series of a function $f(u)$ around 0 converges at $u=1$
 then the $n \rightarrow \infty$ limit of the coefficient of $u^n$ in
 $\frac{f(u)}{1-u}$ is equal to $f(1)$.  \end{proof}

{\bf Remark:} Since the measure $M_{O,u}$ arises for the orthogonal
groups (Theorem 2 of \cite{F2}), Theorem \ref{asymptotics} is a
precise sense in which $gg^{\tau}$ for $g$ random in $GL(n,q)$ is a
hybrid of orthogonal and symplectic groups. Note also that there is a
minor misstatement in Theorem 2 of \cite{F2}: the partitions
$\lambda_{z-1}$ and $\lambda_{z+1}$ do indeed asymptotically both have
the distribution $M_{O,1}$ and are asymptotically independent of
partitions corresponding to other polynomials, but they are not
asymptotically independent of each other.

\section{Bilinear Forms and Conjugacy Classes} \label{bilinearforms}

 Throughout this section the field $F$ of definition of
 $G=GL(n,F)$ and $G^+=\langle G,\tau \rangle$ is arbitrary.

We note that
the conjugacy classes in the coset of $G \tau$ in $G^+$ are in
bijection with the orbits of $G$ acting on $G$ by
$a \circ g = a g a'$.  We will say two matrices are congruent
if they are in the same $\GL$ orbit under this action.

This is a classical problem studied by many (see the bibliography
-- in particular, see Gow \cite{gow1}, \cite{gow2}, \cite{gow3}).
We can view $g$ as defining a nondegenerate bilinear form on
$V:=F^n$ (via $(u,v)=u'gv$). Then the orbits of $G$ are precisely
the congruence classes of nondegenerate bilinear forms. We let
$O_g=O_g(V)$ denote the stabilizer of $g$ in this action (this is
precisely $C_G(g\tau)$, the elements in $G$ which commute with $g
\tau$).

We make some remarks about this problem and relate it to the results in
the earlier part of the article.   It is convenient to use both points
of view (the group theoretic and linear algebra).

We first make some simple observations.  Note that $(g\tau)^2 = gg^{\tr}$.

\begin{lemma} \label{block} Let $g \in G$ and set $h =gg^{\tr}$.
\begin{enumerate}
\item  $(xgx')(xgx')^{\tr} = xhx^{-1}$;
\item   $g^{-1}hg = h^{\tau}$.
\end{enumerate}
\end{lemma}

So we can replace $g$ by a congruent element (and so change $h$ by
conjugation) so that $h$ has a nice form -- for example, we can we
may assume that $h$ is in block diagonal form with the diagonal
blocks corresponding to the primary decomposition for $h$ (i.e.
the characteristic polynomial of the $i$th block is a power of the
irreducible polynomial $f_i$ of degree $d_i$). Let $V_i$ denote
the corresponding subspace.

Moreover, since $h$ is conjugate to $h^{\tr}$ and so to $h^{-1}$
(any element in $G$ is conjugate to its transpose), we may assume
that either the set of roots of $f_i$ is closed under inverses (we
say $f_i$ is a self dual polynomial) or that the roots of $f_j$
are the inverses of the roots of $f_i$.

Lemma \ref{block} implies that conjugating $h$ by $g$ sends $h$ to
$h^{\tau}$ which is still in block diagonal form. Since $h^{\tau}$
is conjugate to $h^{-1}$, $g$ must send a block so its inverse
block. In other words, $g$ preserves each block corresponding to a
self inverse $f_i$ and interchanges the blocks corresponding to
the paired blocks.  This implies:

\begin{lemma} $V$ is an orthogonal direct sum of spaces where the characteristic
polynomial of $h$ is either a power of a self inverse irreducible polynomial
or is a power of $f_1f_2$ where $f_1$ is irreducible with a root $\alpha $ and
$f_2$ is the minimal polynomial of $\alpha^{-1}$.  Moreover, this orthogonal
decomposition is unique (up to order).
\end{lemma}

So $V$ is an orthogonal direct sum of the blocks (or paired
blocks) and there is no loss in assuming that either there is a
single block or precisely two paired blocks.  We now assume that
is the case and consider these summands.

We deal with this last case.  Then $n=2m$. Replacing $g$ by an
equivalent element, allows us to assume that $h =\diag(B,B^{\tr})$
where the characteristic polynomial of $B$ is a power of one of
the two irreducible factors (and the other factor corresponds to
$B^{\tau}$).  This is essentially the primary decomposition for
$h$ -- specifying the relationship between the two blocks. Since
$g$ interchanges the two blocks, we see that all solutions to
$gg^{\tr}=h$ are of the form:
$$ g = \begin{pmatrix} 0 & A \\  D^{\tr} & 0  \end{pmatrix}, $$
where $AD=B$ and $A$ (so also $D$) commutes with $B$.

In particular, a straightforward computation shows that all
solutions to $gg^{\tr}=h$ are in the same $C_G(h)$ orbit.  Thus,
we may replace $g$ and assume that $A=B$ and $D=I$.

Once we have made this simplification, we see that  since $O_g$
centralizes $h$, $O_g = \{ \diag(C,C^{\tr})| BC=CB \}$. Let
$$ J= \begin{pmatrix} 0 & I \\  -I & 0  \end{pmatrix}.$$
Then $h$ preserves the alternating form defined by
$J$ and $O_g = C_{\Sp(J)}(h) \cong C_{GL_m(F)}(B)$.  In particular, we see
that (over an algebraically closed field)  $O_g$
is connected and has dimension at least $m$ (the rank
of $\Sp(J)$ and every centralizer has dimension at least $m$).
Indeed for a generic $g$, we
see that $O_g$ is a maximal torus in this symplectic group.

Note also that all nondegenerate alternating forms preserved by
$h$ are in a single $C_G(h)$ orbit (for $B$ and $B^{\tau}$  have
no fixed points on $\wedge^2$) and so any such form must be of the
form
$$\begin{pmatrix} 0 & X' \\  -X & 0  \end{pmatrix},$$
and so in the orbit of $C_G(h)$. Moreover, the argument above
showed that $C_S(h)$ was independent of the choice of the
symplectic group containing $h$.

Next consider the case where $f$ is irreducible of degree $n=2m$
and the roots of $f$ are closed under inverses.   By passing to a
(finite) Galois extension of $F$, we can reduce to the previous
case (if $F$ is separably closed, every irreducible polynomial has
a single root and so since we are excluding $\pm 1$ as possible
roots, we have split $f$ appropriately). By a general descent
argument (or arguing as in \cite{gow3}) since $O_g$ is a product
of $GL's$ modulo its unipotent radical, we see that in this case,
being equivalent over an extension field implies equivalence over
the original field (in particular over a finite field, we just
apply Lang's theorem). We can always choose our alternating form
to be defined over $F$ (if $F$ is infinite, use a density argument
or just take an $F$-basis for the fixed points on the exterior
square and a generic linear combination over $F$ will be
nondegenerate;  if $F$ is finite, we saw that over the algebraic
closure, these form a single $C_G(h)$ orbit with a stabilizer
being $C_S(h) \cong C_{\GL(m)}(h)$). So as in the previous case,
for any symplectic group $S$ containing $h$ (defined over $F$),
$O_g = C_S(h)$.

We record these observations (see also \cite{gow3} and
\cite{waterhouse}).

\begin{theorem}  Suppose that $h:=(g \tau)^2$ has minimal polynomial
relatively prime to $z^2-1$. Then $n=2m$ is even.  There exist
nondegenerate $h$-stable alternating forms defined over $F$.
Moreover for any such form with corresponding symplectic group
$S$, we have:
\begin{enumerate}
\item $O_g = C_S(h)$ is connected (over the algebraic closure).
\item If $xx^{\tr}$ is conjugate to $h$, then $x$ and $g$
are congruent.
\item If $g,x \in G$ and $g$ and $x$ are congruent over
$\bar F$, then they are congruent over $F$.
\item  Any two solutions $xx^{\tr} = h$ are congruent via an element
of $C_G(h)$;
\item $g$, $g'$ and $g^{-1}$ are all congruent.
\end{enumerate}
\end{theorem}

Thus, if $g \in GL(n,F)$ with $h$ as above, we see that there is a
bijection between similarity classes of such $h$ and classes of
bilinear forms. Since the first situation is well understood, so
is the second.

Using the previous result, we can give another proof of Theorem
\ref{Hallcor} for such elements.

\begin{cor}  Let $F$ be a finite field and
$h \in G=\GL(n,F)$ with the minimal polynomial of
$h$ relatively prime to $z^2-1$.  Assume that $h$ and $h^{-1}$
are conjugate in $\GL(n,F)$. In particular, $n$ is even.
\begin{enumerate}
\item The number of solutions of
$gg^{\tau}=h$ is $|C_G(h):O_g|$;
\item The probability that a random $g \in G$ satisfies
$gg^{\tau}$ is conjugate to $h$ is equal to the probability
that a random element of $\Sp(n,F)$ is conjugate to $h$.
\item The probability that a random $g \in G$ satisfies
$gg^{\tau}$ is regular semisimple is the probability that
a random element of $\Sp(n,F)$ is regular semisimple.
\item The probability that a random $g \in G$ satisfies
$gg^{\tau}$ does not have eigenvalue $\pm 1$ is the probability that
a random element of $\Sp(n,F)$ does not have $\pm 1$ as an eigenvalue.
\end{enumerate}
\end{cor}

\begin{proof}  As we noted above, over the algebraic closure of
$F$, the solutions to $gg^{\tau}=h$ are a single $C_G(h)$ orbit with
connected stabilizer $O_g$.  Over a finite field, we can just apply
Lang's theorem to conclude the same is true over $F$.  This proves the
first statement.

Thus, the number of  solutions $gg^{\tau}$ conjugate to $h$ is
$|G:O_g|$ and the probability that a random $g$ satisfies this is
$1/|O_g|= 1/|C_S(h)|$ as claimed.

Summing over regular semisimple $h \in S$ yields the third
statement.

Summing over all such $h \in S$ yields the final statement.
\end{proof}

If $n > 1$ is odd, we obtain a similar result. Let $F$ be a finite
field of cardinality $q$. Assume that $h$ is conjugate to its
inverse and has determinant $1$. Then $1$ is an eigenvalue for
$h$.  We may assume that $h = \diag(h_0,1)$ where $h_0$ is as in
the previous result. So if $gg^{\tau}=h$, then $g = \diag(g_0,1)$
and $O_g = O_{g_0} \times \mu_2$ where $\mu_2$ is the group of
second roots of unity in $F$. Moreover, we note that we could
replace the symplectic group in the previous theorem by an
orthogonal group (with essentially no change in proof).  There is
an issue of what form of the orthogonal group.  However, in the
case we are over a finite field with $n$ odd this is not a
problem.  Thus, we see that $O_g = C_{O_n(F)}(h)$.

We also see that the number of solutions to the equation to
$xx^{\tau}=h$ is exactly $q-1$ times the number of solutions to
$yy^{\tau}=h_0$. Denote this last number by $N(h_0)$. Thus, the
number of solutions to $xx^{\tau}$ being a conjugate of $h$ is
$|G:C_G(h)|(q-1)N(h_0)$.

Thus, the probability that a random $x$ is a solution to
$xx^{\tau}$ conjugate to $h$ is $(q-1)N(h_0)/|C_G(h)|=
N(h_0)/|C_{\GL(n-1,q)}(h_0)|$ is the probability that a random
$y \in \GL(n-1,q)$ is a solution to $yy^{\tau}=h_0$.  Summarizing
we have:

\begin{cor}  Let $F$ be a finite field and
$h \in G=\GL(n,F)$ with $h = \diag(h_0,1)$ with
the characteristic polynomial
$f$ of $h_0$
relatively prime to $z^2-1$.  Assume that $h$ and $h^{-1}$
are conjugate in $\GL(n,F)$. In particular, $n=2m + 1$ is odd.
\begin{enumerate}
\item $O_g = C_{O_n(F)}(h)$ for some orthogonal group $O_n(F)$
containing $h$;
\item The number of solutions of
$gg^{\tau}=h$ is $|C_G(h):O_g||\mu_2(F)|$;
\item The probability that a random $g \in G$ satisfies
$gg^{\tau}$ is conjugate to $h$ is equal to the probability
that a random element of $\Sp(n-1,F)$ is conjugate to $h_0$.
\item The probability that a random $g \in G$ satisfies
$gg^{\tau}$ is regular semisimple is the probability that
a random element of $\Sp(n-1,F)$ is regular semisimple.
\item The probability that a random $g \in G$ satisfies
$gg^{\tau}$ has characteristic polynomial $(z-1)w(z)$
where $w(z)$ is relatively prime to $z^2-1$ is the probability
that a random element of $\Sp(n-1,F)$ has no eigenvalue $\pm 1$.
\end{enumerate}
\end{cor}

Note that many of the results stated above do not hold for all nondegenerate
bilinear forms.  Any two symmetric invertible matrices are congruent over
an algebraically closed field of characteristic not $2$
but this is not the case for fields which
are not  quadratically closed.

We now point out some consequences over an algebraically closed
field.

\begin{lemma}  Let $F$ be an algebraically closed field.
Set $G:= \GL(n,F)$ and $E:=
\{ g \in G: gg^{\tr} \text{ has distinct eigenvalues} \}$.
Suppose that $g \in E$.
\begin{enumerate}
\item
$E$ is a nonempty open subset of $G$;
\item Either $n$ is even and
$O_g$ is a torus of dimension $n/2$ or $n=2m+1$ and
$O_g = T \times \mu_2$ where $T$ is a torus of dimension $m$
and $\mu_2$ is the groups of second roots of unity in $F$.
\item  There exists an involution $x \in G$ such that the
$1$-eigenspace of $x$ has dimension $\lfloor n/2 \rfloor$ with $xgx'=g'$.
\item Any solution $x$ to $xgx'=g'$ is an involution.  Moreover,
either $x$ or $-x$ has fixed space of dimension $\lfloor n/2 \rfloor$.
\end{enumerate}
\end{lemma}

\begin{proof}  $E$ is clearly open.  Arguing as we have before, we
see that $V$ is an orthogonal direct sum of subspaces of dimension
$2$ plus a $1$-dimensional summand if $n$ is odd (each two dimensional
summand is the span of $h$-eigenvectors corresponding to
inverse eigenvalues --  the $1$-dimensional summand
is the fixed space of $h$).

To show it is nonempty thus reduces to the   $2 \times 2$ case, where it is clear.
Again, the computation of $O_g$ reduces to checking the cases where
$n \le 2$.  If $n=1$, then $O_g = \mu_2$ and if $n=2$, $O_g$ is a
$1$-dimensional torus.   Similarly, the last statement reduces to the
case of $n \le 2$.
\end{proof}

\begin{cor}  $\dim O_A \ge \lfloor{n/2} \rfloor$ for any $A$.
\end{cor}

\begin{proof}  This is true on an open dense subvariety by the previous
result.  It is now standard to see the result holds on the closure
(consider the subvariety $(g,A) \in G \times M_n(F)$ where $gAg'=A$;
then $O_A$ is the fiber over $A$ of the projection onto the second fiber
and has dimension $[n/2]$ on an open dense subset of $M_n(F)$, whence
on all of $M_n(F)$).
\end{proof}

Next consider the case that the characteristic polynomial of
$h$ is $(z + 1)^n$ (assume that $F$ does not have characteristic $2$).
We can view $G^+$ as a subgroup of $GL(2n,F)$ and can consider the Jordan
decomposition of $g\tau$.  Then $g\tau = u g_1\tau = g_1 \tau u$ where
$u$ is unipotent and $(g_1 \tau)^2 = -1$; i.e. $g_1$ is a skew symmetric
matrix.  By passing to a congruent element, we may assume that $g_1 = J$, some fixed skew
symmetric matrix (and again $n$ is even -- or use the fact that
$h$ has determinant $1$). The fact that $u$ commutes with $J \tau$ is equivalent
to the fact that $u \in \Sp_J$.  So we see that such elements correspond
to unipotent conjugacy classes in the symplectic group (and conversely
since we are not in characteristic $2$, every such element is the square
of $g \tau$ for some $g$).

So in this case, there is a bijection between unipotent conjugacy classes
in the symplectic group and equivalence classes of such forms.

Finally, consider the case that the characteristic polynomial of
$h$ is $(z - 1)^n$.   If $F$ does not have characteristic $2$, we
argue precisely as above and we see that we may take $g \tau = u g_1 \tau =
g_1 \tau u$ where $g_1$ is a symmetric matrix and $u$  is a unipotent
element in the orthogonal group corresponding  to  the symmetric bilinear
form $g_1$.   In  particular, if $F$ is finite, there are two choices
for the class of $g_1$ and then  the unipotent classes in each of
the corresponding orthogonal  groups.

If $F$ has characteristic $2$, then $g \tau$ has order a power of $2$
and so $g \tau$ is contained in a maximal unipotent subgroup of $G^+$.
All  such are conjugate (essentially by Sylow's theorem or its analog
for linear groups) and so we see that we may assume that $g \in U$,
a maximal unipotent subgroup of $G$ and $\tau$ normalizes $U$.

These observations will show that:
\begin{lemma}  Assume that $F$ is a finite field.
  Write $n=2m + \delta$ with $\delta$ either
$0$ or $1$.
If $h:=(g \tau)^2$ has characteristic polynomial $(z \pm 1)^{n}$, then
either $O_g$ has order divisible by $q^m$ or $q$ is odd, $n=2$
and $|O_g| > q^m$.
\end{lemma}

\begin{proof} If $n=1$, there is nothing to prove. So assume that $n \ge 2$.

  First consider the case where
$F$ has odd characteristic.
Then $O_g = C_{H}(h^2)$ where $H$ is a symplectic or
orthogonal group
containing the unipotent element $h^2$.  If $H$ is
a symplectic group, then $\delta=0$ and $H$ has rank
$m$.  This result is well known for semisimple groups
(either by inspection of the classes or by counting
fixed points on unipotent subgroups -- cf \cite{FG1}).
If $H$ is an orthogonal group, the same argument applies
unless $n=2$ (the semisimple rank of $H$ is $m$).
If $n=2$ and $H$ is an orthogonal group, then $O_g$
is $O_2^{\epsilon}(q)$ and so has order greater than $q$
(but not divisible by $q$).

Now suppose that $F$ has characteristic $2$ (a variant of this
approach would work in characteristic not $2$ as well). We may assume
that $g \tau $ normalizes the standard unipotent subgroup $U$
and $h \in U$ (because any two Sylow $2$-subgroups of $G^+$
are conjugate).

We claim that $C_U(g \tau)$ has order divisible by $q^m$.
It suffices to show that after passing to the algebraic closure
that $\dim C_{\bar U}(g\tau) \ge m$ (here ${\bar U}$ is the maximal
unipotent subgroup containing $U$ over the algebraic closure).
This is because the fixed points of the $q$-Frobenius map on any
connected unipotent group of dimension $m$ has $q^m$ elements.
Let $V=g \tau {\bar U}$.  This is a connected variety.  Let
$s:V \rightarrow U$ be the squaring map.  It suffices
to show that $\dim C_{\bar U} ( g \tau u) \ge m$ for $u$
in an nonempty open subset of ${\bar U}$.
If $n=2$, then the Sylow $2$-subgroup of $G^+$ is abelian
and the result is clear.  So assume that $n > 2$.
Suppose that $n > 2$ is odd.  Let $R$ be the set of
regular unipotent elements in ${\bar U}$.

  By Theorem \ref{Hallcor}, it follows that $s^{-1}(R)$ is nonempty
  and so is a nonempty open subvariety.  If $g\tau u$ is in this
  set, we claim that $O_g=C_{\GL}(g\tau u) \le {\bar U}$.  The centralizer
  of $(g \tau u)^2$ is contained in $T{\bar U}$ where $T$ is the
  group of scalars.  Note that $C_T(g \tau u)=1$, whence the claim.
  If $n$ is even, we replace $R$ by $R_1$, the set of elements
  in ${\bar U}$ that correspond to a partition of shape $(n-1,1)$.
  Again by Theorem \ref{Hallcor}, the set of $g \tau u$ whose
  square is in $R_1$  is a nonempty open subvariety.  For such an
  element, $C_{\GL}( g \tau u)  \le T{\bar U}$ where $T$ is
  a $2$-dimensional torus.  Moreover, $T$ has eigenspaces of dimension
  $1$ and $n-1$, whence  $g \tau u$ acts as inversion on $T{\bar U}/{\bar U}$
  and so $C_{\GL}( g \tau u) \le \bar{U}$.
  So in either case, we have shown that there is a nonempty open
  subset of $g\tau {\bar U}$ whose centralizer is contained in $\bar U$.
  We have already observed that any centralizer in $\GL$ has dimension
  at least $m$, whence the result. \end{proof}

These results immediately yield a completely different proof of the
lower bound in Theorem \ref{centralizersize}:

\begin{cor} \label{cent} Let $F$ be a finite field.
Let $g \in \GL(2m + \delta,F)$ with $\delta =0$ or $1$.
The minimum size of $O_g$ is at least the smallest centralizer size of an
element in $|\Sp(2m,F)|$, and hence at least
$(1-\frac{1}{q^2}-\frac{1}{q^4})^2 q^{m} \left(\frac{1-1/q}{4elog_{q}(2m)} \right)^{1/2}$.
\end{cor}
\begin{proof}  Set $h=(g \tau)^2$. We split $V$ as an orthogonal sum
of $V_i, 1 \le i \le 3$ -- where $h$ is unipotent on $V_1$, $h$
is $-u$ with $u$ unipotent on $V_2$ and the minimal polynomial
of $h$ is relatively prime to $(z^2-1)$ on $V_3$.
On $V_2 \oplus V_3$, the previous result implies that
$|O_g|$ is at least as big as the centralizer of some element
in $Sp(n_2 + n_3,q)$ (note that $n_i:=\dim V_i$ is even for $i \ne 1$).

On $V_1$,  we see that $|O_g| \ge q^{\lfloor {n_1/2} \rfloor}$.
If $n_1$ is even, take an element in $Sp(n,q)$ (note $n$ is even)
that is regular unipotent of size $n_1 + n_2$ and $h$ on $V_3$
and we see that this centralizer is no bigger than $|O_g|$.
If $n_2$ is odd, take the element as above in $\Sp(n-1,q)$
and conclude the same result.

    The lower bound on centralizer sizes of elements in $Sp(2m,q)$
    appears in \cite{FG3}. \end{proof}

\section{Number of Conjugacy Classes} \label{numberclasses}

 This section gives upper bounds for the number of $G^+(n,q)$
 conjugacy classes in the coset $GL(n,q) \tau$ and also treats
 a variation (which we use in \cite{FG2}) for $SL(n,q)$. Let us
 make some preliminary remarks about $G^+(n,q)$ to show that
 our bound has substance. It is well-known that $GL(n,q)$ has
 at most $q^n$ conjugacy classes, and so it follows that the
 number of conjugacy classes in $G^+(n,q)$ is at most
 $2q^n$. In fact we shall see that the number of classes in the
 coset $GL(n,q) \tau$ is at most $28 q^{\lfloor n/2
 \rfloor}$. Throughout this section, we will let $k(GL(n,q)
 \tau)$ denote the number of $G^+(n,q)$ conjugacy classes in
 the coset $GL(n,q) \tau$.

 In fact Gow \cite{gow3} derived (in the language of bilinear forms)
 generating functions for the number of $G^+(n,q)$ conjugacy classes in the
 coset $G(n,q) \tau$. See Waterhouse \cite{waterhouse} for a different
 proof of Proposition \ref{genfunc}. They did not however, give explicit upper bounds.

\begin{prop} (\cite{gow3}) \label{genfunc} Let $g(t)$ be the generating function
\[ g(t) = 1+\sum_{n \geq 1} t^n k(GL(n,q) \tau).\]
Let $f=1$ if the characteristic is even and $f=2$ if the characteristic is odd.
Then $g(t)=\prod_{i \geq 1} \frac{(1+t^i)^f}{1-qt^{2i}}.$ \end{prop}

\begin{lemma} \label{NPlem} For $q \geq \sqrt{2}$,
\[ 1-\frac{1}{q}-\frac{1}{q^2}+\frac{1}{q^5}+\frac{1}{q^7}-\frac{1}{q^{12}}-\frac{1}{q^{15}}
< \prod_{i \geq 1 } (1-\frac{1}{q^i}) < 1-1/q.\]
\end{lemma}

\begin{proof} This is proved along the same lines as Lemma 3.5 of \cite{NP}.
Namely the upper bound is obvious and the lower bound follows from Euler's
pentagonal number theorem (exposed in \cite{A}), which states that
\begin{eqnarray*}  \prod_{i \geq 1}(1-1/q^i) & = & 1+\sum_{n \geq 1}
(-1)^n(q^{-n(3n-1)/2}+q^{-n(3n+1)/2})\\
& = & 1-1/q-1/q^2+1/q^5+1/q^7-1/q^{12}-1/q^{15}\\
& & + 1/q^{22}+1/q^{26}
-\cdots. \end{eqnarray*}
Since $q \geq \sqrt{2}$, one has that the sum
consisting of powers of $q$ higher than $1/q^{26}$ has magnitude less
than $1/q^{22}$.
\end{proof}

\begin{lemma} (\cite{MR}) \label{MR} The coefficient of $t^n$ in
$\prod_{i \geq 1} \frac{1-q^i}{1-tq^i}$ is at most $q^n$.
\end{lemma}

We remark that the generating function in Lemma \ref{MR}
is the generating function for the number of conjugacy classes in $GL(n,q)$.
Theorem \ref{classdouble} gives upper bounds on $k(GL(n,q)\tau)$.

\begin{theorem} \label{classdouble}
\begin{enumerate}
\item $k(GL(n,q) \tau) \leq 28q^{\lfloor n/2 \rfloor}$ if $q$ even.
\item $k(GL(n,q) \tau) \leq 23q^{\lfloor n/2 \rfloor}$ if $q$ odd.
\end{enumerate}
\end{theorem}

\begin{proof} Let us first consider $k(GL(2n,q) \tau)$ in the case that
the characteristic is even. By Proposition \cite{gow3}, for the first part
we seek the coefficient of $t^{2n}$ in $\prod_{i} \frac{1+t^i}{1-qt^{2i}}$.
Writing this generating function as \[\prod_{i} \frac{1-t^{2i}}{1-qt^{2i}}
\prod_i \frac{1 }{1-t^{i}}, \] one sees from Lemma \ref{MR} and the fact
that all coefficients in $\prod_i \frac{1+t^i}{1-t^{2i}}$ are positive that
the sought coefficient is at most
\begin{eqnarray*} & & q^n \sum_{m \geq 0} \frac{1}{q^m} Coeff.
\ t^{2m} \ in \ \prod_i \frac{1}{1-t^{i}}\\
& \leq & q^n \prod_i \frac{1}{1-1/q^{i/2}}. \end{eqnarray*}
The
quantity $\prod_i \frac{1}{1-1/q^{i/2}}$ is maximized (among legal
$q$) for $q=2$. Then Lemma \ref{NPlem} implies that $k(GL(2n,q) \tau)
\leq 28 q^n$. The case $k(GL(2n+1,q) \tau)$ with even characteristic
is similar.

 For the second part, let us examine the case $k(GL(2n,q) \tau)$
 where the characteristic is odd. Writing the generating function as
 \[ \prod_{i \geq 1} \frac{1-t^{2i}}{1-qt^{2i}} \prod_{i \geq 1}
 \frac{1+t^{i}}{1-t^i} \]
 the same argument shows that the sought
coefficient is at most
\begin{eqnarray*} & & q^n \sum_{m \geq 0} \frac{1}{q^m} Coeff.
\ t^{2m} \ in \ \prod_i \frac{1+t^{i}}{1-t^{i}}\\ & \leq & q^n \prod_i
\frac{1+1/q^{i/2}}{1-1/q^{i/2}}. \end{eqnarray*}
The quantity $\prod_i
\frac{1+1/q^{i/2}}{1-1/q^{i/2}}$ is maximized (among legal $q$) for $q=3$.
Rewriting this as $\prod_i
\frac{1-1/q^{i}}{(1-1/q^{i/2})^2}$ and applying Lemma \ref{NPlem}
(once to upper bound the numerator and once to lower bound the denominator)
establishes the upper bound of $23q^n$. The case $k(GL(2n+1,q) \tau)$ is similar.
\end{proof}

 Next we determine the fixed $q$, large $n$ asymptotics of $k(GL(n,q) \tau)$.

\begin{lemma} \label{bign} (Darboux \cite{O}) Suppose that $f(u)$
is analytic for $|u|<r,r>0$ and has a finite number of simple poles on
$|u|=r$. Letting $w_j$ denote the poles, and $g_j(u)$ be such that
$f(u)=\frac{g_j(u)}{1-u/w_j}$ and $g_j(u)$ is analytic near $w_j$,
then as $n \rightarrow \infty$, the difference between the coefficient of
$u^n$ in $f(u)$ and $\sum_j \frac{g_j(w_j)}{w_j^n}$ goes to 0.
\end{lemma}

\begin{prop} \label{asy} Suppose that $q$ is fixed. Let $f=1$ if the characteristic is
even and $f=2$ if the characteristic is odd.
\begin{enumerate}
\item
\[ lim_{n \rightarrow \infty} \frac{k(GL(2n,q) \tau)}{q^n} =
\frac{1}{2} \frac{\prod_{i \geq 1} (1+\frac{1}{q^{i/2}})^f +
\prod_{i \geq 1} (1+\frac{(-1)^i}{q^{i/2}})^f}{\prod_{i \geq 1} (1-\frac{1}{q^i})} \]
\item
\[ lim_{n \rightarrow \infty} \frac{k(GL(2n+1,q) \tau)}{q^{n}} =
\frac{q^{.5}}{2} \frac{\prod_{i \geq 1} (1+\frac{1}{q^{i/2}})^f -
\prod_{i \geq 1} (1+\frac{(-1)^i}{q^{i/2}})^f}{\prod_{i \geq 1} (1-\frac{1}{q^i})} \]
\end{enumerate}
\end{prop}

\begin{proof} Let us prove the first part, the second part being similar.
>From Proposition \ref{genfunc}, $\frac{k(GL(2n,q) \tau)}{q^n}$
is the coefficient of $t^{2n}$ in
\[ g(\frac{t}{\sqrt{q}}) = \frac{1}{1-t^2} \prod_{i \geq 1} \frac{(1+t^i/q^{i/2})^f}{1-t^{2(i+1)}/q^i}.\]
The result now follows from Lemma \ref{bign}.
\end{proof}

 Although we do not need it, we include the following proposition for
 completeness.

\begin{prop} \label{total}
\begin{enumerate}
\item $k(G^+(n,q)) = \frac{1}{2} k(GL(n,q)) + \frac{3}{2} k(GL(n,q)\tau)$.
\item For $q$ fixed and $n$ big, $k(G^+(n,q))$ is asymptotic is $\frac{q^n}{2}$.
\end{enumerate}
\end{prop}

\begin{proof} As explained in \cite{gow1},
a real conjugacy class of $GL(n,q)$ remains a conjugacy class in
$G^+(n,q)$ and an inverse pair of non-real conjugacy classes of $GL(n,q)$
merges into a single conjugacy class in $G^+(n,q)$. Thus the elements of
$GL(n,q)$ account for $k(GL(n,q))/2$ plus one half the number of real
conjugacy classes of $GL(n,q)$. The first assertion now follows from
\cite{gow3}, which shows that $k(GL(n,q) \tau)$ is the number of real
conjugacy classes of $GL(n,q)$.

    The second assertion follows from the first assertion, together with Proposition
    \ref{asy} and the fact from \cite{FG3} that $k(GL(n,q))$ is asymptotic to $q^n$ for $q$ fixed.
\end{proof}

 Next we treat $<SL(n,q),\tau>$ classes rather than $<GL(n,q),\tau>$
 classes.

We consider the orbits of $\SL(n,F)$ on $M_n(F)$ under the action
$A \rightarrow gAg'$.  Clearly, $\SL(n,F)$ preserves determinant.
Moreover, given any $A$, there is certainly is $B$ in the
$\GL$-orbit of $A$ with $\det(B) = b^2 \det(A)$ for any $b \in F$.
So we restrict our attention to matrices with a fixed determinant
(and all that matters is the square class of the determinant). So
consider $A \in \GL(n,F)$ with $\det(A) = d$ nonzero.

\begin{lemma}  The set of matrices of determinant $d \ne 0$ that
are $GL$ congruent to $A$ is a single $\SL$-orbit if
$O_A$ contains an element of determinant $-1$ and splits into
two orbits otherwise (which are in bijection).
\end{lemma}

We note that for $n$ even a generic $A$ (i.e. $AA^{\tau}$ having distinct
eigenvalues) satisfies $O_A \le \SL(n,F)$ ($O_A$ is a torus contained
in a symplectic group).  So in this case the  $\GL$ orbit
of $A$ intersect the set of matrices with $\det = d$ splits into
two orbits for $\SL$.  If $n$ is odd, $-I \in O_A$ for every $A$.

Thus, over a finite field:

\begin{lemma}  Let $F$ be a finite field and $n$ a positive integer.
Fix $g \in \GL(n,F)$.
\begin{enumerate}
\item If $n$ is odd, then the number of $\SL(n,F)$ conjugacy classes
in the coset $\SL(n,F) g \tau$ is the number of $\GL(n,F)$ conjugacy
classes in the coset $\GL(n,F)\tau$.
\item If $n$ is even, then the number of $\SL(n,F)$ conjugacy
classes in the coset $\SL(n,F) g \tau$ is less than twice
number of $\GL(n,F)$ conjugacy classes in the coset $\GL(n,F)\tau$.
\end{enumerate}
\end{lemma}

We also see that the $\langle \SL(n,q),\tau \rangle $ centralizer of $g \tau$
is either equal to the $\langle \GL(n,q),\tau \rangle$ centralizer
or has index $2$ in it (with equality in characteristic $2$ or generically
when $n$ is even and always nonequality
if $qn$ is odd).
 So the bounds from Section \ref{mincentsize} are applicable. Also, the smallest
centralizer size does not change when $n$ is even.

\section{Some Examples with Derangements} \label{derangements}

In the series of papers beginning with \cite{FG1}, the authors
have verified Shalev's conjecture that the proportion of
derangements in a simple group is bounded away from $0$ (by an
absolute constant).  This immediately implies the same result for
$G^{+}(n,q)$.  However, we give some examples to show that if we
restrict our attention to the coset containing $\tau$ (which is
relevant to images of rational points on curves over finite fields),
this need not be the case. For a full
treatment, see the forthcoming paper \cite{FG2}. Our purpose here
is to give some examples which are instant corollaries of results
in earlier sections of this paper.

Set $G=\GL(n,q)$ and $G^+ = \langle G, \tau \rangle$.  All our actions
are projective actions and so we are really working in the quotient.
In particular, for $n=2$, $\tau$ is an inner automorphism on
$\PGL(2,q)$.  We thus assume that $n > 2$. \\

\noindent
{\bf Example 1 }   Suppose that $n > 2$ is even. Let $\Omega$ be the
set of $1$-dimensional spaces of alternating nondegenerate forms
over $F_q$.   This is a single $G$ orbit and is acted on by
$G^{+}$ with stabilizer $H$ the normalizer of $GSp(n,q)$.  Let
$S=Sp(n,q) < H$.

  Suppose that
$h:=(x \tau)^2$ is a regular semisimple element.  We have seen
that $C_G(x \tau)$ is conjugate to a maximal torus of $S$.

Now suppose that $q$ is even.  It is straightforward to see that
we can solve $(y\tau)^2 =h$ with $y \tau$ normalizing $S$ (because
every semisimple element has odd order).   Then $x \tau$ and $y
\tau$ are conjugate, whence $x\tau$ has a fixed point on $\Omega$.
By Corollary \ref{regss} the probability that $(x \tau)^2$ is
regular semisimple is equal to the probability that an element of
$Sp(n,q)$ is regular semisimple, and hence approximately $1 - 2/q$
(\cite{GL},\cite{FNP}) when $q$  is big.   So for large $q$, we
see that the proportion of derangements in the coset of $\tau$ is
at most approximately $2/q$ and goes to $0$ as $q \rightarrow
\infty$ (independently of $n$).

We can say a bit more.  Let $F$ be the algebraic closure of $F_q$.
In characteristic $2$, the non semisimple regular elements in
$\Sp(n,F)$ form a subvariety of codimension $1$ with $2$
components (one consisting of elements that commute with a long
root element and one with a short root element).  One computes
that the generic element $h$ of the component consisting of elements
that commute with a short root element has no eigenvalue $1$. By
our earlier results, such elements are of the form $gg^{\tau}$.  On
the other hand, there is no solution $xx^{\tau}=h$ with $x \in GSp(n,F)$
(for $\Sp(n,F)$ is the centralizer of $J \tau$ for some skew symmetric
matrix $J$ and so $ (s(J \tau))^2 = h$ implies that $h$ is a square
in $GSp(n,F)$ which is easily seen not to be the case).  Thus,
the elements in this component generically are derangements and so
we see that the proportion of derangements is at least $O(1/q)$.

    For a lower bound with $q$ fixed, note that if $xx^{\tau}$ is
    regular semisimple (i.e. square-free characteristic polynomial)
    except for having a non self-conjugate pair $\{\phi, \bar{\phi}\}$
    of degree 1 polynomials which each have Jordan type consisting of
    a single part of size 2, then $x \tau$ is a derangement. Now
    suppose that $q>2$ even is fixed. Then by Theorem
    \ref{asymptotics}, the $n \rightarrow \infty$ proportion of such
    elements is \[ \frac{q-2}{2} \frac{1}{q^2(1-1/q)}
    \frac{rs_{Sp}(n,\infty)}{1+\frac{1}{q-1}}\] (the factor
    $\frac{q-2}{2}$ counts the number of possible pairs $\{\phi,
    \bar{\phi}\}$). Here $rs_{Sp}(n,\infty)$ is the fixed $q$ large $n$
    limiting proportion of regular semisimple elements in $Sp(n,q)$,
    proved in \cite{FNP} to lie between
    $\frac{(q-1)^2(q^2+2q+2)}{q^2(q+1)^2}$ and $\frac{q-1}{q+1}$ when
    $q$ is even. Hence the limiting proportion of such $x \tau$ is
    bounded away from $0$ for small $q$ and large $n$; moreover the
    lower bound is roughly $\frac{1}{2q}$ for $q$ not too small. For a
    lower bound for $q=2$, note that if $xx^{\tau}$ has a z-1
    component of dimension 4 and has Jordan structure 3,1, then it
    cannot be in a symplectic group. By Theorem \ref{asymptotics}, the
    fixed $q$ large $n$ limiting probability that
    $\lambda_{z-1}(gg^{\tau})$ has size 4 with parts of size 1 and 3
    is $\frac{1}{q^2} \prod_{i=1}^{\infty} (1-\frac{1}{q^{2i-1}})$.


If $q$ is odd, then a similar analysis shows that if $gg^{\tau}=h$ is
a regular semisimple element and is a nonsquare, then $g \tau$ is a
derangement in this action.  Since every maximal torus of
$\Sp(2n,q)$ has even order, at most $1/2$ (and typically much less)
of the elements of a maximal torus are squares.  If $q$ is large, then
almost all elements are regular semisimple and at least close to
$1/2$ of them are nonsquares.  Thus, the limiting proportion of derangements
in the coset of $\tau$ is at least $1/2$ as $q \rightarrow \infty$.

    For fixed $q>3$, again note that if $xx^{\tau}$ is regular
semisimple except for having a non self-conjugate pair $\{\phi,
\bar{\phi}\}$ of degree 1 polynomials which each have Jordan type
consisting of a single part of size 2, then $x \tau$ is a derangement.
Then by Theorem \ref{asymptotics}, the $n \rightarrow \infty$
proportion of such elements is \[ \frac{q-3}{2} \frac{1}{q^2(1-1/q)}
\frac{rs_{Sp}(n,\infty)}{1+\frac{1}{q-1}}\] (the factor
$\frac{q-3}{2}$ counts the number of possible pairs $\{\phi,
\bar{\phi}\}$). Here $rs_{Sp}(n,\infty)$ is the fixed $q$ large $n$
limiting proportion of regular semisimple elements in $Sp(n,q)$,
proved in \cite{FNP} to lie between
$1-\frac{3}{q}+\frac{5}{q^2}-\frac{10}{q^3}$ and
$1-\frac{3}{q}+\frac{5}{q^2}-\frac{6}{q^3}$ when $q$ is odd. The case
$q=3$ can be treated exactly as the case $q=2$. Hence the limiting
proportion of such $x \tau$ is bounded away from $0$ for fixed $q$ and
large $n$.

    To summarize, we see that either the proportion of
derangements in the coset $G \tau$ goes to $0$ (as $2/q$ does) in the
case $q$ is even or is bounded away from $0$ if $q$ is odd. Actually
when $q$ is odd we proved uniform boundedness away from $0$ for all
but finitely many $(n,q)$; however uniform boundedness for all $n,q$
now follows by a result of Jordan that any transitive permutation
group acting on a set of size $n>1$ has a derangement.\\

\noindent
{\bf Example 2 } The next example shows that for $n > 2$ odd, there is an action
with few derangements in the coset of $\tau$.

Let $\mathcal{E}$ denote the set of $g$ such that $h:=gg^{\tau}$
has characteristic polynomial
$(z-1)^{\epsilon}w(z)$ where $w(z)$ is prime to $z^2-1$ and
$\epsilon = 0$ if $n$ is even and $1$ if $n$ is odd. If $n$ is
odd, $h$ is conjugate to $\diag(h_0,1)$.  Let $S=S_h$ be a
symplectic group containing $h$ (or $h_0$ if $n$ is odd).

\begin{lemma}  Suppose that $n > 2$ is odd.  Let $G^+$ act
on the set  $\Gamma$ of complementary point-hyperplane pairs.  Then $g \in
\mathcal{E}$ implies that $g \tau$ has a fixed point on $\Gamma$.
\end{lemma}


\begin{proof}  Set $h=gg^{\tau}$.
We can embed $h \in H:=GL(1,q) \times GL(n-1,q)$ and then we see
that $xx^{\tau}=h$ has a solution with $x \in H$. Thus, $x \tau$
normalizes $H$ and so $x \tau$ has a fixed point, whence $g \tau$
does (as $g \tau$ is conjugate to $ \lambda x \tau$ for some
$\lambda \in F_q$).
\end{proof}


In particular, if $gg^{\tau}$ is semisimple regular, this implies that
$g$ has a fixed point.  The proportion of such $g \tau$ with square
not being regular semisimple is roughly $2/q$ (for $q$ even) and $3/q$
for $q$ odd (see Lemma \ref{Hallcor} and \cite{GL} or \cite{FNP}).  So
the proportion of derangements in the coset $G \tau$ goes to zero as
$q \rightarrow \infty$
 (independently of $n$). Also note from Theorem \ref{asymptotics} that the
fixed $q$, $n \rightarrow \infty$ proportion of elements in
$\mathcal{E}$ is $\prod_{j \geq 1} (1-\frac{1}{q^{2j-1}})^f$ where
$f=2$ if the characteristic is odd and $f=1$ if the characteristic is
even. For $q$ not too small this is roughly $1-\frac{f}{q}$.

    For a lower bound in the case of fixed $q$, $n \rightarrow
    \infty$, note that if $gg^{\tau}$ has $\lambda_{z-1}$ being one
    part of size 3 and the characteristic polynomial of
    $gg^{\tau}$ has no degree 1 factors other than $z-1$, it must
    be a derangement on the set of complementary point-hyperplane
    pairs. By Theorem \ref{asymptotics}, the proportion of such
    elements is \[ \frac{1}{q(1-1/q^2)} \prod_{j \geq 1}
    (1-\frac{1}{q^{2j-1}})^f (1-\frac{1}{q^j})^{(q-1-f)/2}\] where
    $f=2$ if the characteristic is odd and $f=1$ if the
    characteristic is even. Since $q \geq 2$, this is at least
    $\frac{c}{q}$ where $c$ is a constant which is easy to make
    explicit.\\

\noindent
{\bf Example 3 } We now consider other actions on pairs of subspaces.
For convenience we assume that $n > 2$ is even (a similar analysis
suffices for $n$ odd other than the case above).  Fix $k < n-k$. Let
$\Omega_k$ be the set of complementary pairs of subspaces of
dimension $k$ and $n-k$.  Let $\Gamma_k$ be the set of flags of
type $k,n-k$ (i.e pairs of subspaces $U_1 \subset U_2$ where $\dim
U_1=k$ and $\dim U_2=n-k$).  Note that $G^+$ acts on both of these
sets and that $G$ acts transitively.


\begin{lemma}  Assume that $n$ is even
and  $g \in \mathcal{E}$ and set
$h = gg^{\tau}$.  Assume also that $h$ is semisimple.
\begin{enumerate}
\item $g \tau$ has a fixed point on $\Omega_k$ if and only
if $h$ fixes a nondegenerate $k$-dimensional subspace (with
respect to the alternating form defining $S$).
\item
$g \tau$ has a fixed point on $\Gamma_k$ if and only
if $h$ fixes a totally singular $k$-dimensional subspace (with
respect to the alternating form defining $S$).
\end{enumerate}
\end{lemma}

\begin{proof} Consider the various actions and let $H$ be
the stabilizer of a point in one of these representations.

If $xx^{\tau}=h$, then $x \tau$ and $g\tau$ are conjugate,
so it suffices to show that such an $x$ exists with $x$
fixing a point precisely when $h$ satisfies the conditions.

Consider $\Omega_k$.  The stabilizer of a point is
$\langle H, \tau \rangle$ where $H=GL(k) \times GL(n-k)$.
If $x \in H$ and $xx^{\tau}$ is conjugate to $H$, then
$h$ must be real on both $k$ and $n-k$ dimensional space.
This implies that $h$ fixes a nondegenerate $k$ dimensional
subspace (with respect to any $h$-invariant alternating form).

Conversely, if $h$ does fix a nondegenerate subspace of dimension
$k$, there is a conjugate of $h$ in $H$ real in both
$GL(k)$ and $GL(n-k)$, whence the result.

The proof of the second assertion is similar (note also that
a semisimple element of $\Sp$ fixes a totally singular
$k$-dimensional subspace implies that it fixes a nonsingular
subspace of dimension $2k$).
\end{proof}

It is proved in \cite{FG4} that for all but finitely many symplectic
groups, the proportion of elements which are regular semisimple and
derangements on totally singular or nondegenerate $k$-spaces is
bounded away from 0 by an explicit absolute constant. From Jordan's
theorem mentioned in Example 1, it follows that the proportion of
derangements in this example is uniformly bounded away from $0$.


\begin{thebibliography}{AAA}




\bibitem[A]{A} Andrews, G., The theory of partitions, Cambridge
University Press, Cambridge, 1984.

\bibitem[BKS]{BKS} Bannai, E., Kawanaka, N., and Song, S., The
character table of the Hecke algebra $H(GL_{2n}(F_q),Sp_{2n}(F_q))$,
{\it J. Algebra}, {\bf 129} (1990), 320-366.

\bibitem[B]{B} Boston, N., Dabrowski, W., Foguel, T.,  et al,
The proportion of fixed-point-free elements of a transitive
permutation group, {\it Comm. Algebra}  {\bf 21} (1993),
3259--3275.





\bibitem[F1]{F1} Fulman, J., Random matrix theory over finite fields,
{\it Bull. Amer. Math. Soc.} {\bf 39} (2002), 51-85.




\bibitem[F2]{F2} Fulman, J., A probabilistic approach to conjugacy classes in
the finite symplectic and orthogonal groups, {\it J. Algebra} {\bf 234} (2000), 207-224.




\bibitem[F3]{F3} Fulman, J., Cycle indices in the finite classical groups,
{\it J. Group Theory} {\bf 2} (1999), 251-289.




\bibitem[FG1]{FG1} Fulman, J. and Guralnick, R., Derangements in simple and primitive groups,
to appear in {\it Durham 2001 Conference on Groups, Combinatorics, and Geometry}.




\bibitem[FG2]{FG2} Fulman, J. and Guralnick, R., Derangements in
almost simple groups, in preparation.




\bibitem[FG3]{FG3} Fulman, J. and Guralnick, R.,
Derangements in classical groups for non subspace actions, preprint.

\bibitem[FG4]{FG4} Fulman, J. and Guralnick, R.,
Derangements in subspace actions of finite classical groups, preprint.


\bibitem[FNP]{FNP} Fulman, J., Neumann, P. and Praeger, C.,
A generating function approach to matrix enumeration in the finite classical
groups, preprint.




\bibitem [Ga]{gabriel}
Gabriel, P., Appendix: Degenerate bilinear forms. {\it J. Algebra} {\bf 31}
(1974), 67-72.




\bibitem[Go1]{gow1} Gow, R., Properties of the characters of the
finite general linear group related to the transpose-inverse involution,
{\it Proc. London Math. Soc.} {\bf 47} (1983), 493-506.




\bibitem[Go2]{gow2} Gow, R., The equivalence of an invertible matrix to
its transpose. {\it Linear and Multilinear Algebra} {\bf 8} (1979/80),
329-336.




\bibitem[Go3]{gow3} Gow, R. The number of equivalence classes of
nondegenerate bilinear and sesquilinear forms over a finite field.
{\it Linear Algebra Appl.} {\bf 41} (1981), 175-181.


\bibitem[GL]{GL} Guralnick, R. and L\"ubeck, F.,
On $p$-singular elements in Chevalley groups in characteristic $p$.
Groups and computation, III (Columbus, OH, 1999), 169--182, Ohio State Univ.
Math. Res. Inst. Publ., 8, de Gruyter, Berlin, 2001.

\bibitem[He]{H} Herstein, I.N., Topics in algebra, Second edition,
Xerox College Publishing, 1975.




\bibitem[HZ]{HZ} Howlett, R. and Zworestine, C.,
On Klyachko's model for the representations of finite general linear groups,
in {\it Representations and quantizations}, 229-245, China High. Educ. Press, Beijing, 2000.




\bibitem[IS]{IS} Inglis, N. and Saxl, J., An explicit model for the complex
representations of the finite general linear groups,
{\it Arch. Math. (Basel)} {\bf 57} (1991), 424-431.




\bibitem[Ka]{Ka} Kawanaka, N., Symmetric spaces over finite fields, Frobenius-Schur indices,
and symmetric function identities, in {\it Physics and combinatorics 1999},
World Scientific, Singapore, 2001.




\bibitem[Kl]{K} Klyachko, A., Models for complex representations of groups
$GL(n,q)$, {\it Math. USSR-Sb.} {\bf 48} (1984), 365-379.




\bibitem[M]{M} Macdonald, I., Symmetric functions and Hall polynomials,
Second edition. Clarendon Press, Oxford, 1995.




\bibitem[MR]{MR} Maslen, D. and Rockmore, D.,
Separation of variables and the computation of Fourier transforms on
finite groups, I., {\it J. Amer. Math. Soc.} {\bf 10} (1997), 169-214.




\bibitem[NP]{NP} Neumann, P. and Praeger, C., Cyclic matrices over finite fields,
{\it J. London Math. Soc. (2)} {\bf 52} (1995), 263-284.




\bibitem[Od]{O} Odlyzko, A.M., Asymptotic enumeration methods, Chapter 22 in
{\it Handbook of Combinatorics, Volume 2}. MIT Press and Elsevier, 1995.




\bibitem[Ok]{Ok} Okounkov, A., Symmetric functions and random partitions,
in {\it Symmetric functions 2001: Surveys of developments and perspectives},
Kluwer Academic Publishers, 2002.




\bibitem [Ri]{riehm}  Riehm, C., The equivalence of bilinear forms,
{\it J. Algebra} {\bf 31} (1974), 45--66.




\bibitem[Se]{S} Serre, J.-P.  Linear representations of finite groups,
Springer-Verlag, New York, 1977.






\bibitem[Wal]{W} Wall, G.E., On conjugacy classes in the unitary, symplectic,
and orthogonal groups, {\it J. Austr. Math. Soc.} {\bf 3} (1963), 1-63.




\bibitem [Wat]{waterhouse} Waterhouse, W., The number of congruence
classes in $M_n(\mathbf{F}_q)$. {\it Finite Fields Appl.} {bf 1} (1995),
57-63.




\end{thebibliography}
\end{document}